\documentclass[1p]{elsarticle}
\usepackage{amsmath} 
\numberwithin{figure}{section}
\numberwithin{table}{section}
\numberwithin{equation}{section}
\usepackage[hidelinks]{hyperref}
\usepackage[hyphenbreaks]{breakurl}
\usepackage{eurosym}
\DeclareRobustCommand{\officialeuro}{%
  \ifmmode\expandafter\text\fi
  {\fontencoding{U}\fontfamily{eurosym}\selectfont e}}
\usepackage{amssymb}
\usepackage[printonlyused]{acronym}
\usepackage{longtable}
\usepackage{caption}
\usepackage[]{units}
\usepackage{mydefs}
\usepackage{arydshln}

\journal{European Journal of Operational Research}

\begin{document}

\begin{frontmatter}

\title{Multiplicity of equilibria in conjectural variations models of natural gas markets}

\author[IED,IfA]{Tobias Baltensperger\corref{correspondingauthor}}
\cortext[correspondingauthor]{Corresponding author. Tel.: +41 44 632 44 73;}
\ead{t.baltensperger@usys.ethz.ch}

\author[ZHAW]{Rudolf M. F\"{u}chslin}
\author[IED]{Pius Kr\"{u}tli}
\author[IfA]{John Lygeros}

\address[IED]{Institute for Environmental Decisions (IED), ETH Z\"{u}rich, 8092 Z\"{u}rich, Switzerland}
\address[IfA]{Automatic Control Laboratory (IFA), ETH Z\"{u}rich, 8092 Z\"{u}rich, Switzerland}
\address[ZHAW]{Institute of Applied Mathematics and Physics (IAMP), ZHAW Z\"{u}rich University of Applied Sciences, 8401 Winterthur, Switzerland}

\begin{abstract}

Spatial partial equilibrium models incorporating conjectural variations are widely used to analyze the development of oligopolistic multi-agent markets, such as international energy and raw material markets. Although this model type can produce multiple equilibria under commonly used assumptions, to the best of our knowledge, the consequences for the interpretation of the model results have not yet been explored in detail.
To this end, we derive a linear complementarity model for the gas market and discuss under which assumptions on the model structure a component of the solution is unique. In particular, we find that the gas flow between a trader and a consumer is unique whenever the trader is modeled to exert market power in the consumer's market. 
We demonstrate our findings by computing the extreme points of the polyhedral solution space and show that erroneous conclusions could be drawn whenever only one (arbitrary) point in the solution space is picked for interpretation. Furthermore, we discuss whether economically meaningful parameter value changes exist which would enforce uniqueness in all components of the solution.

\end{abstract}


\end{frontmatter}

\setcounter{MaxMatrixCols}{20}
\section{Introduction} \label{sec:P1_Introduction}
Spatial partial equilibrium models including conjectural variations (SPE-CV models) are a popular tool to investigate oligopolistic markets, such as fossil fuel or electricity markets \citep{Neuhoff2005,Diaz2010,Kaminski2011}. These markets are characterized by a high capital intensity as a result of high exploration costs, long lead times, large distances between producing and consuming regions, or expensive infrastructure requirements. This creates high barriers to entry and limits the competition, which allows major traders to exert market power over consumers.

In SPE-CV models, multiple traders are distinguished to represent the impact of market power exertion. As a consequence, one of the model outcomes are the sales of the individual traders to consumers, which is of interest in practice, because it indicates which two parties are likely to trade in a given situation. Unfortunately, this representation may also introduce redundancies to the model, mathematically corresponding to the existence of multiple solutions for a given set of parameters, and economically to the existence of multiple market equilibria in a single situation. This is problematic when it comes to the interpretation of the model results, since one could not determine which, of the possibly many, equilibria proposed by the model become reality.

In the presented study we investigate SPE-CV models for gas markets. Several variants have been proposed in the past, including the special cases of Cournot and Bertrand competition. The most prominent examples are the dynamic GASTALE model by \citet{Lise2008b}, and the \ac{WGM} by \citet{Egging2008}: both models assume international gas traders to exert market power over consumers, and compete perfectly for infrastructure services such as transmission or storage capacity. The dynamic GASTALE model \citep{Lise2008,DeJoode2010} and the \ac{WGM} \citep{Zhuang2008,Gabriel2009,Egging2010,Egging2010a,Gabriel2012,Egging2013,Huppmann2013a,Moryadee2014} were expanded successively and used for various case studies.

However, in most articles, the multiplicity of solutions is not discussed. Exceptions include \citet{Gabriel2005a}, who analyzed a predecessor of the \ac{WGM} and derived uniqueness criteria for the intermediate and wholesale prices, and \citet{Egging2010a}, who presented an example in which non-unique service prices can arise in the \ac{WGM}. To the best of our knowledge, the remaining variables, particularly the volume flows, have not yet received any attention. The contribution of this article is to analyze the multiplicity of all variables by exploiting the information hidden in the mathematical structure of the problem. For the non-unique variables, we discuss whether uniqueness is relevant for interpretation in the first place, and investigate how uniqueness could be achieved based on changes in the model's parameter values.

The remainder of this article is structured as follows. In the next section, we introduce a SPE-CV model for the gas market. In Section \ref{sec:P1_ExistanceAndUniqueness}, we substantially reduce the number of equations by expressing (infrastructure) service prices and volume flows explicitly. Service cost functions are affine or quadratic in the volume flows, allowing the model to be represented as a linear complementarity problem. Subsequently, we carry out the mathematical analysis concerning the multiplicity of solutions. In Section \ref{sec:P1_NumericalExample}, we demonstrate our findings via a numerical example. In Section \ref{sec:P1_Discussion}, we discuss our theoretical and numerical results, and in Section \ref{sec:P1_Conclusions}, we summarize and conclude our work.

\section{Model description} \label{sec:P1_ModelDescription}
\subsection{General setting and notation} \label{sec:P1_GeneralSetting}
The structure of our SPE-CV model comprises nodes and arcs, and distinguishes multiple time periods within a year. The nodes represent countries or regions, in which gas is produced, stored, and sold to consumers. The arcs represent the gas flow capacities from one region to another. The differentiation of time periods allows for a distinction of multiple levels of demand. The geographical and temporal granularity and range of the model can be adjusted to fit the requirements of a specific case; in the numerical example presented in Section \ref{sec:P1_NumericalExample}, the world-wide gas trade over a year is represented by $50$ nodes, $291$ arcs, and $2$ time periods (summer and winter half-year). 

While SPE-CV models are certainly able to capture the main properties of the gas market such as limited production, restricted transport capabilities, large seasonal demand fluctuations, and market power exerting traders, they do not provide any information on the status of the gas market within a node, implicitly assuming that the gas is distributed well and prices are uniform. If the region mapping to a node is chosen sufficiently small, this weakness can be compensated to some extent. Since publicly available data for some of the parameters, particularly production capacities and wholesale prices, is mostly limited to country-level, a finer granularity is rarely aimed for with SPE-CV models. Other types of models, such as the TIGER model \citep{Lochner2011a,Lochner2011,Lochner2012,Dieckhoner2013}, or the model introduced by \citet{Carvalho2014}, circumvent this data constraint by assuming a fixed and therefore price-independent demand; this eliminates the need for accurate wholesale price data and hence allows for a higher model resolution. 

\begin{figure}[!htb]
\centering
\includegraphics[width=\textwidth]{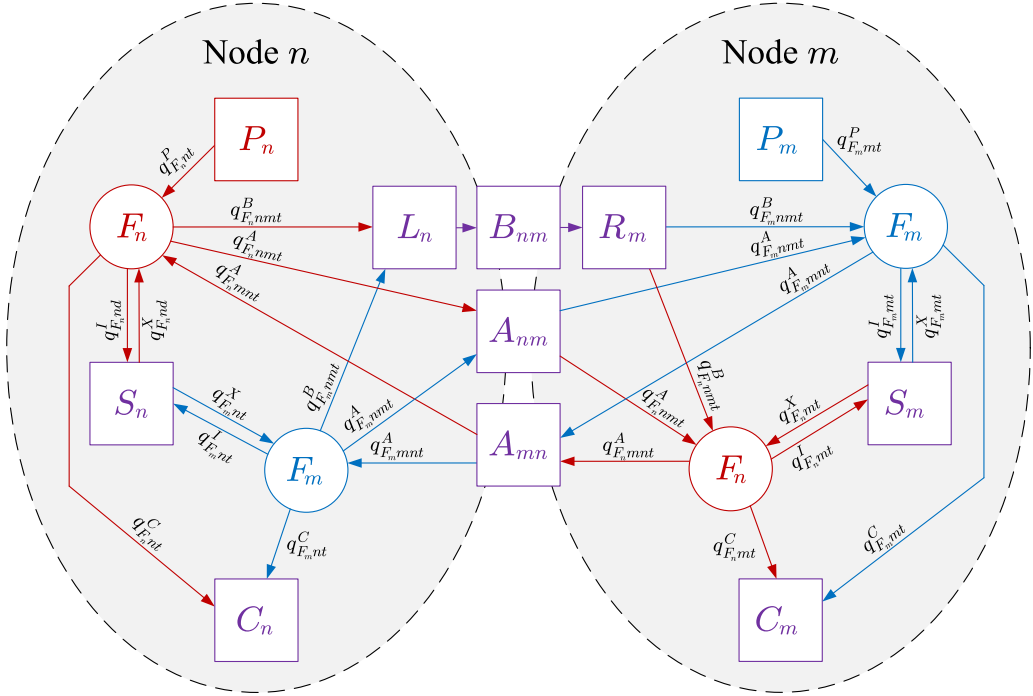}
\caption{Gas market model with two nodes. $P$: producer. $A$: pipeline operator. $L$: liquefaction plant operator. $B$: LNG shipment. $R$: regasification plant operator. $S$: storage operator. $C$: consumer. $F_n$: trader associated with producer $P_n$. $F_m$: trader associated with producer $P_m$. $q^{P}_{fnt}$: quantity delivered from producer to trader $f$ at node $n$ in time period $t$. $q^{A}_{fnmt}$: pipeline transportation of trader $f$ via arc $nm$ in period $t$. $q^{B}_{fnmt}$: LNG shipment of trader $f$ via arc $nm$ in period $t$. $q^{I}_{fnt}$: storage injection by trader $f$ at node $n$ in period $t$. $q^{X}_{fnt}$: storage extraction by trader $f$ at node $n$ in period $t$. $q^{C}_{fnt}$: sales by trader $f$ to consumer in node $n$ in period $t$. The traders $F_n$ and $F_m$, there decision variables, and their corresponding producers $P_n$ and $P_m$ are colored red ($n$) and blue ($m$), respectively. Service providers, except producers, and flows between them are marked purple, as well as the consumers, since all traders trade with them.}
\label{fig:P1_P1_Model_full}
\end{figure}

The mathematical formulation is inspired by other models in the field, particularly the deterministic version of the \ac{WGM} proposed by \citet{Egging2010a}, except that we refrain from modeling investments over time. In comparison to the \ac{WGM}, we aim for a more compact representation of the overall model to facilitate the mathematical analysis, and therefore choose a more general notation for the individual components of the model to emphasize their similarity. The notation is introduced as we formulate the model, and summarized in Tables \ref{tab:P1_ServiceProviders}-\ref{tab:P1_Functions}. An exemplary model with two interconnected nodes is depicted in Figure \ref{fig:P1_P1_Model_full}.

The nodes $n, m \in \mathcal{N}, n \neq m$ represent countries or regions, and the time periods are denoted $t \in \mathcal{T}$. In each node $n$, a consumer $C_n$, a gas producing company $P_n$, a storage operator $S_n$, a liquefaction plant operator $L_n$, and a regasification plant operator $R_n$ may be located. The storage operator $S_n$ provides two services: injection of gas ($I_n$), and extraction of gas ($X_n$). A transmission system operator $A_{nm}$ manages the flows in the pipelines from $n$ to $m$, and a shipping company $B_{nm}$ transports \ac{LNG} from $n$ to $m$. To simplify notation in the remainder of this paper, we introduce $Z_z$ as a placeholder for a producer or a infrastructure service $Z$ in a node/arc $z$, thus $Z \in \{P,I,X,L,R,A,B\}$, and $z \in \mathcal{Z}$, where $\mathcal{Z}$ is the corresponding subset of nodes/arcs in which a service\footnote{Although producers are not an infrastructure service, we henceforth include them in our notion of \textit{services}.} of type $Z$ is active. Note that each service $Z_z$ represents all companies of type $Z$ in a node/arc $z$, even if in reality multiple firms carry out this activity in the geographical region of $z$. Similarly, $C_n$ represents the aggregated demand of all consumers located in node $n$.

In most of the large companies, gas is produced in one division, and then internally transferred to the trading-arm of the company, which deals with the (international) trade and distribution. Each of these companies' trading-arm is modeled by a trader $F_n$, which receives gas only from the producer $P_n$ located in its home-node. Each $F_n$ spreads out from his home-node through the network, thereby competes with the other traders for service capacities (except production capacities), and finally sells gas in the consumer markets, where it again competes with the other traders for market shares. While we assume all services to be price-takers, the traders exert market power in the consumer markets, which is modeled based on a conjectural variations approach. 

\subsection{Mathematical derivation}
The traders and service providers have perfect and complete information about the market. Therefore, the decision processes of the traders and the service providers can be modeled as (deterministic) optimization problems. In the following, we introduce the optimization problem of the traders (Section \ref{sec:P1_OptTrader}), the optimization problem of the service providers (Section \ref{sec:P1_OptProblemOfSysServ}), and the market clearing conditions (Section \ref{sec:P1_MarketCC}).

\subsubsection{Optimization problem of the traders} \label{sec:P1_OptTrader}
The goal of each trader $f$ is to maximize profit $\Pi^F_{f}$, which is equal to total revenue $\TotRev^F_{f}$ from sales to the consumers, minus the total cost $\TC^F_{f}$ arising from using services:
\begin{equation}  
\Pi^F_{f} = \TotRev^F_{f} - \TC^F_{f}.
\end{equation}
$f \in \mathcal{F}$ represents an arbitrary trader in the set of all traders (whereas $F_n$ is the specific trader whose home base is in node $n$). The total revenue of trader $f$ reads
\begin{equation} \label{eqn:P1_TraderRevenue}
\TotRev^F_f = \sum\limits_{t \in \mathcal{T}} \sum\limits_{n \in \mathcal{C}(f)} \left[\left(1-\theta_{fnt}^C \right) \lambda^C_{nt} + \theta_{fnt}^C \Lambda^{C}_{nt}(\cdot) \right]  q^{C}_{fnt}. 
\end{equation} 
The decision variables are all the $q^{C}_{fnt}$ which represent the total amount of gas sold to the consumers in node $n$ and period $t$. The quantity is multiplied by the corresponding wholesale price, and the product is summed over all periods $\mathcal{T}$ and over the set of nodes $\mathcal{C}(f)$ in which the trader $f$ and a consumer are present. 
The wholesale price comprises the weighted sum of two quantities: the price on the wholesale market $\lambda^{C}_{nt}$, which is determined by the market clearing conditions and is exogenous to the trader; and the inverse demand function $\Lambda^{C}_{nt}$, a function of the total demand $\sum \limits_{f \in \mathcal{F}(n)} q^C_{fnt}$, on which the trader has influence via its decision variable $q^C_{fnt}$. The weight $\theta^C_{fnt} \in [0,1]$ is referred to as the behavioral parameter of market power \citep[Chapter 12]{Tremblay2012} of trader $f$ over consumers located in node $n$ in period $t$: if $\theta^C_{fnt}=1$ for all traders $f$, the market $n$ in period $t$ is characterized by a Cournot equilibrium, whereas $\theta^C_{fnt}=0$ for all traders $f$ indicates a competitive equilibrium. Conceptually following the setting used by \citet{Egging2010} for the \ac{WGM}, we also allow values $0<\theta^C_{fnt}<1$, representing intermediate stages of the traders' conjectures about the other traders' behaviors\footnote{There are different propositions how values of $\theta^C_{fnt}$ different from 0 and 1 should be interpreted, see for example \citet[Chapter 3]{Perloff2007} for a discussion. In this work, we follow the interpretation of \citet[Chapter 12]{Tremblay2012} and think of $\theta^C_{fnt}$ as "toughness of competition", where decreasing $\theta^C_{fnt}$ indicates increasing competition.}. However, we do not allow $\theta^C_{fnt}>1$ and hence exclude cartelization of multiple traders.

The inverse demand function in node $n$ and period $t$ is defined as
\begin{equation} \label{eqn:P1_LambdaC} 
\Lambda^{C}_{nt}(\cdot) := \INT^C_{nt} + \SLP^C_{nt} \sum \limits_{f \in \mathcal{F}(n)} q^{C}_{fnt}.
\end{equation}
$\INT^C_{nt}$ represents the intercept of the inverse demand function, which corresponds to the maximum willingness to pay of consumers at node $n$ in period $t$, $\SLP^C_{nt}$ is the slope and is assumed strictly negative. These parameters are derived from the reference willingness to pay $\WTP^C_{nt}$, reference demand $\DMD^C_{nt}$, and price elasticity $\eta_{nt}$, which in turn is calculated from the price elasticities of the residential \& commercial, industrial and power sectors and their reference shares $\SHR^C_{nd,\cdot}$ in consumption. 
\begin{align}
\INT^C_{nt}=& (1-\frac{1}{\eta^C_{nt}})\cdot \WTP^C_{nt} \geq 0\\
\SLP^C_{nt}=&\frac{\WTP^C_{nt}}{\DMD^C_{nt} \cdot \eta^C_{nt}} < 0\\
\eta^C_{nt} =& \eta^C_{nd,res} \cdot \SHR^C_{nd,res} + \eta^C_{nd,ind} \cdot \SHR^C_{nd,ind} + \eta^C_{nd,el} \cdot \SHR^C_{nd,el} < 0 \label{eqn:P1_eta_nt}
\end{align}
All parameters can vary in $n$ and $t$ if the respective data is available. Following previous works of \citet{Lise2008} (dynamic GASTALE model) and \citet{Egging2010} (\ac{WGM}), the inverse demand function is chosen affine; for a market analysis including more general structures of inverse demand functions we refer to \citet{Abolhassani2014}.

For each trader $f$, the total cost arising from service utilization can be written as 
\begin{equation} \label{eqn:P1_TraderCost}
\TC^F_f = \sum \limits_{\mathcal{Z} \in \{\mathcal{P},\mathcal{A},\mathcal{R},\mathcal{B},\mathcal{L},\mathcal{I},\mathcal{X}\}} \sum\limits_{t \in \mathcal{T}} \sum\limits_{z \in \mathcal{Z}(f)} \lambda^Z_{zt} q^{Z}_{fzt},
\end{equation} 
where the $q^{Z}_{fzt}$ are the decision variables of the traders and represent the volume flows between trader $f$ and service provider $Z_z$, $\lambda^Z_{zt}$ is the market price for using the service $Z_z$ (exogenous to the trader), and $\mathcal{Z}(f)$ is the set of nodes/arcs at which the trader $f$ and a service provider of type $Z$ are both present. Note that 
\begin{align*}
q^{L}_{fnt} &\equiv \sum \limits_{m \in \mathcal{B}(n)} \frac{ q^{B}_{fnmt}}{{\LOSS^L_{n}}} \text{, and} \\
q^{R}_{fmt} &\equiv \sum \limits_{n \in \mathcal{B}(m)} \LOSS^B_{nm} q^{B}_{fnmt},
\end{align*}
where $q^{L}_{fnt}$, $q^{B}_{fnmt}$, and $q^{R}_{fmt}$ are the flows at arrival at the service providers, and $\LOSS^L_{n}$ and $\LOSS^B_{nm}$ are the loss factors from liquefaction and shipment; $q^{L}_{fnt}$ and $q^{R}_{fmt}$ will therefore be substituted for the remainder of this work.

Each trader's profit maximization is subject to several volume balances to ensure that the solution is reasonable: per trader $f$, all inflows have to match the outflows in each node $n$ and each period $t$, and the yearly inflows into a storage unit have to be equal to the outflows in each node $n$. 
\begin{align} 
\begin{split}
h^{F1}_{fnt}=& q^{P}_{fnt} + q^{X}_{fnt} + \sum \limits_{m \in \mathcal{A}(n)} \LOSS^A_{mn} q_{fmnt}^{A} + \sum \limits_{m \in \mathcal{B}(n)} \LOSS^{B}_{mn} \LOSS^{R}_{n} q_{fmnt}^{B} \\
&- q^{I}_{fnt} -  q^{C}_{fnt} - \sum \limits_{m \in \mathcal{A}(n)} q_{fnmt}^{A} - \sum \limits_{m \in \mathcal{B}(n)} \frac{q_{fnmt}^{B}}{\LOSS^L_n} = 0\quad \left(\phi_{fnt}^{N} \right) \quad \forall n,t  \end{split} \\
h^{F2}_{fn}=& \sum \limits_{t \in \mathcal{T}} \LOSS^I_n q^{I}_{fnt} - \sum \limits_{t \in \mathcal{T}} q^{X}_{fnt} = 0 \quad \left(\phi^{S}_{fn} \right)\quad \forall n  \label{eqn:P1_VolumeBalances_2}
\end{align}
Some incoming and outgoing flows are multiplied and divided, respectively, by $\LOSS^Z_z \in \left(0,1 \right]$, to account for incurring transaction losses. $\phi_{fnt}^{N}$ and $\phi^{S}_{fn}$ are the associated dual variables (Lagrange multipliers).

Finally, all the decision variables $q^C_{fnt}$ and $q^Z_{fzt}$ can be upper and lower bounded due to contracts, certain policies, sanctions, etc. We state these equations in the following, but for simplicity of notation will neglect them for the remainder of the paper. 
\begin{align} 
q^C_{fnt} \leq& \UPP^C_{fnt} \quad \left(\overline{\xi}^C_{fnt} \right) \quad \forall f,n,t \label{eqn:P1_UpperCLimits}\\
q^C_{fnt} \geq& \LOW^C_{fnt} \quad \left(\underline{\xi}^C_{fnt} \right) \quad \forall f,n,t \label{eqn:P1_LowerCLimits}\\
q^Z_{fzt} \leq& \UPP^Z_{fzt} \quad \left(\overline{\xi}^Z_{fzt} \right) \quad \forall f,z,t \label{eqn:P1_UpperZLimits}\\
q^Z_{fzt} \geq& \LOW^Z_{fzt} \quad \left(\underline{\xi}^Z_{fzt} \right) \quad \forall f,z,t \label{eqn:P1_LowerZLimits}
\end{align}
\subsubsection{Optimization problem of the service providers} \label{sec:P1_OptProblemOfSysServ}
Each service provider $Z_z$ determines the total volume contracted to traders $s^Z_{zt}$ by maximizing profits $\Pi^Z_z$ over all time periods $t$:
\begin{align} \label{eqn:P1_SysServOpt}
\max\limits_{s^Z_{zt}\geq 0} \quad \Pi^Z_{z} &= \sum\limits_{t \in \mathcal{T}} \left(\lambda^Z_{zt} s^Z_{zt}- c^Z_{zt}(s^Z_{zt}) \right)\\
\text{s.t.} \quad		g^{Z1}_{zt} &= \overline{\CAP}^Z_{zt} - s^Z_{zt} \geq 0 \quad \left(\alpha^Z_{zt} \right) \quad \forall t \\
								  g^{Z2}_{z} &= \overline{\CAP}^{ZT}_{z} - \sum \limits_{t \in \mathcal{T}} s^Z_{zt} \geq 0 \quad \left(\alpha^{ZT}_{z} \right). \label{eqn:P1_g2}
\end{align}
$c^Z_{zt}(s^Z_{zt})$ are the costs arising for the service provider $Z_z$ by contracting $s^Z_{zt}$ and are assumed to be convex and twice differentiable, $\overline{\CAP}^{Z}_{zt}$ and $\overline{\CAP}^{ZT}_{z}$ are the service providers' capacity limits per time period $t$ and per year, and $\alpha^Z_{zt}$ and $\alpha^{ZT}_{z}$ are the dual variables associated with the capacity constraints. As stated above, $\lambda^Z_{zt}$ is the market price for using the service and is exogenous to the service provider. 

Note that $g^{X2}_n$ can be used to model the working gas volume of a storage: $\overline{\CAP}^{XT}_{n} - \sum \limits_{t \in \mathcal{T}} s^X_{nt} \geq 0$, where $\overline{\CAP}^{XT}_{n}$ is the maximum storable working gas in node $n$. Limiting the total extractable gas is equal to limiting the total storable gas under the assumption that the injection and extraction periods are separated, and the maximum storage level is reached at most once per year. This is clearly true when modeling up to $\overline{t}=3$ periods per year, but mostly also holds for $\overline{t} > 3$ in reality. Also, the traders have perfect foresight over the entire simulation period, and therefore it is not necessary to specify whether the withdrawal or the injection period is first in the simulation. It suffices to ensure that the storage levels are equal or higher at the end of a storage cycle than before, which is guaranteed in equilibrium by Equations \eqref{eqn:P1_VolumeBalances_2} and \eqref{eqn:P1_g2}. To model the storage more accurately, additional constraints can be introduced in the formulation without affecting the theoretical conclusions, as long as they are affine in the decision variables. 

\subsubsection{Market clearing conditions} \label{sec:P1_MarketCC}
These conditions couple the decision variables of the traders ($q^{Z}_{fzt}$) and the service providers ($s^Z_{zt}$). By enforcing market clearance in every node $n$ and period $t$, the dual variables of these equations, the market prices $\lambda^Z_{zt}$, are determined:
\begin{equation} \label{eqn:P1_MarketClearing}
h^{M1}_{zt} = s^Z_{zt} - \sum \limits_{f \in \mathcal{F}(z)} q^{Z}_{fzt} = 0  \quad \left(\lambda^Z_{zt} \right).
\end{equation}
Furthermore, in all nodes $n$ and periods $t$ we enforce clearing of the wholesale markets:
\begin{equation} \label{eqn:P1_MarketClearingEndConsumer}
h^{M2}_{nt} = \lambda^C_{nt} - \Lambda^C_{nt}(\cdot) = 0 \quad \left(\lambda^C_{nt} \right).
\end{equation}

\subsection{Mixed complementarity problem formulation}
We derive the Karush-Kuhn-Tucker conditions of the optimization problems of the traders and service providers, and complement them with the market clearing conditions. This leads to Problem \eqref{eqn:P1_FullProblem} which takes the format of a mixed non-linear complementarity program. 
\begin{subequations} \label{eqn:P1_FullProblem}
\begin{alignat}{5}
0 &\leq& -\frac{\partial \Pi^F_f}{\partial q^Z_{fzt}} -  \sum \limits_{n^\prime \in \mathcal{N}(f)} \sum \limits_{t^\prime \in \mathcal{T}} \phi^N_{fn^\prime t^\prime} \frac{\partial h^{F1}_{fn^\prime t^\prime}}{\partial q^Z_{fzt}} -  \sum \limits_{n^\prime \in \mathcal{N}(f)} \phi^S_{fn^\prime} \frac{\partial h^{F2}_{fn^\prime}}{\partial q^Z_{fzt}} & \perp q^Z_{fzt} &&\geq 0 \quad \forall f,z,t \\
0 &\leq& -\frac{\partial \Pi^F_f}{\partial q^C_{fnt}} - \sum \limits_{n^\prime \in \mathcal{N}(f)} \sum \limits_{t^\prime \in \mathcal{T}} \phi^N_{fn^\prime t^\prime} \frac{\partial h^{F1}_{fn^\prime t^\prime}}{\partial q^C_{fnt}} -  \sum \limits_{n^\prime \in \mathcal{N}(f)} \phi^S_{fn^\prime} \frac{\partial h^{F2}_{fn^\prime}}{\partial q^C_{fnt}}&\perp q^C_{fnt} &&\geq 0  \quad \forall f,n,t \\
0 &=&  h^{F1}_{fnt} &\perp \phi^N_{fnt} && \left(\text{free} \right) \quad \forall f,n,t \label{eqn:P1_f1}\\
0 &=&  h^{F2}_{fn} &\perp \phi^S_{fn} &&\left(\text{free} \right)  \quad \forall f,n \label{eqn:P1_f2}\\
0 &\leq& -\frac{\partial \Pi^Z_z}{\partial s^Z_{zt}} - \sum \limits_{t^\prime \in \mathcal{T}} \alpha^Z_{zt^\prime} \frac{\partial g^{Z1}_{zt^\prime}}{\partial s^Z_{zt}} -  \alpha^{ZT}_{z} \frac{\partial g^{Z2}_z}{\partial s^Z_{zt}}&\perp s^Z_{zt} &&\geq 0 \label{eqn:P1_g0} \quad \forall z,t \\
0 &\leq& g^{Z1}_{zt} & \perp \alpha^Z_{zt} &&\geq 0  \quad \forall z,t  \\
0 &\leq& g^{Z2}_z & \perp \alpha^{TZ}_z &&\geq 0  \quad \forall z  \\
0 &=& h^{M1}_{zt} &\perp \lambda^Z_{zt} && \left(\text{free} \right)  \quad \forall z,t \label{eqn:P1_h1}\\
0 &=& h^{M2}_{nt} &\perp \lambda^C_{nt} && \left(\text{free} \right)  \quad \forall n,t \label{eqn:P1_h2}
\end{alignat}
\end{subequations}
A solution to Problem \eqref{eqn:P1_FullProblem} is simultaneously an optimal solution to the problems of all the traders and service providers. The solution represents a Nash equilibrium, in which all traders and service providers optimize their volume flows given the other traders' and service providers' decisions, and no party has an incentive to move away from the equilibrium.

\section{Existence and uniqueness of solutions} \label{sec:P1_ExistanceAndUniqueness}
In this section, we address the mathematical characteristics of Model \eqref{eqn:P1_FullProblem}, particularly the existence and uniqueness properties of its solution. Before we start, we restate the model in a more compact way by eliminating superfluous variables, specifying the cost functions, and setting all loss factors $\LOSS^Z_z =1$. 

\subsection{Compact model representation} \label{sec:P1_reformulation}
First, we want to reduce the number of variables and equations by substituting the decision variables of the service providers $s^Z_{zt}$ and the service prices $\lambda^{Z}_{zt}$ and eliminating Equations \eqref{eqn:P1_g0} and \eqref{eqn:P1_h1}. This saves programming time, simplifies maintenance, and potentially reduces calculation times. Since $h^{M1}_{zt}=0$ for all arcs/nodes $z$ and time periods $t$, we can immediately set $s^Z_{zt}:=\sum \limits_{f \in \mathcal{F}(z)} q^{Z}_{fzt}$ for all $z,t$ and eliminate Equation \eqref{eqn:P1_h1}. The wholesale prices $\lambda^{C}_{zt}$ and Equation \eqref{eqn:P1_h2} could be substituted/eliminated in a similar way, however, we refrain from applying these changes to obtain a more favorable mathematical representation of the model for the derivation of our results. 

Eliminating the service prices $\lambda^{Z}_{zt}$ and Equation \eqref{eqn:P1_g0} is somewhat more involved. We consider the complementarity constraint \eqref{eqn:P1_g0} and distinguish two cases: for a specific $z,t$, either
\begin{enumerate}[(i)]
	\item $s^Z_{zt} > 0 \implies \lambda^Z_{zt} = \frac{dc^Z_{zt}(s^Z_{zt})}{ds^Z_{zt}} + \alpha^Z_{zt} + \alpha^{ZT}_{z}$, or 
	\item $s^Z_{zt} = 0 \implies \lambda^Z_{zt} \leq \frac{dc^Z_{zt}(s^Z_{zt})}{ds^Z_{zt}} + \alpha^Z_{zt} + \alpha^{ZT}_{z} = \left. \frac{dc^Z_{zt}(s^Z_{zt})}{ds^Z_{zt}}\right|_{s^Z_{zt}=0}$, the marginal costs for service provider $Z_z$ at zero throughput. $\alpha^Z_{zt}=0$ and $\alpha^{ZT}_{z}=0$ under the assumption that $\overline{\CAP}^Z_{zt}>0$ and $\overline{\CAP}^{ZT}_{z}>0$, which is reasonable, since otherwise the variable $s^Z_{zt}$ would not have been introduced in the first place.
\end{enumerate}
Hence, we set the service price $\lambda^Z_{zt} := \frac{dc^Z_{zt}(s^Z_{zt})}{ds^Z_{zt}} + \alpha^Z_{zt} + \alpha^{ZT}_{z}$ for all $z,t$ regardless of the contracted volume flow $s^Z_{zt}$ being greater or equal to zero. This is an easily justifiable change to the model, since the only expected change in the output is that the service price $\lambda^Z_{zt}$ is also unambiguous for nonoperating services $Z_z$. 

Second, for homogeneity of notation we relax Equalities \eqref{eqn:P1_f1}, \eqref{eqn:P1_f2}, and \eqref{eqn:P1_h2} to $h^{F1}_{fnt} \geq 0$, $h^{F2}_{fn} \geq 0$, and $h^{M2}_{nt} \geq 0$, and simultaneously force non-negativity on the respective dual variables. This reformulation does not change the outcome of the model: we expect $h^{F1}_{fnt}= h^{F2}_{fn}=0$ in equilibrium, because the trader would otherwise be wasting gas, which is never optimal because production costs are assumed positive. Furthermore, in equilibrium we either have $h^{M2}_{nt}=0$ and wholesale prices $\lambda^C_{nt} \geq 0$ (as before), or $h^{M2}_{nt} > 0$ and $\lambda^C_{nt} = 0$ (instead of $h^{M2}_{nt} = 0$ and $\lambda^C_{nt} < 0$). But since it is never optimal for traders to sell gas at $\lambda^C_{nt} \leq 0$ due to production cost being larger than zero, the second situation never arises. In reality, these assumptions may not always hold: it may be cheaper for a supplier to flare superfluous gas or sell gas at negative prices than reducing production. However, in our framework we do not impose any constraints which would lead to such behavior (as for instance ramp constraints), and therefore $h^{F1}_{fnt}\neq 0$, $h^{F2}_{fn} \neq 0$, and $h^{M2}_{nt} \neq 0$ never result in our model, even after relaxation.

Third, we define the cost functions of the service providers $c^Z_{zt}(s^Z_{zt})$ such that sufficiency and necessity of the Karush-Kuhn-Tucker conditions is preserved. Hence, $c^Z_{zt}(s^Z_{zt})$ have to be convex and differentiable. In recent studies \citep{Gabriel2012, Moryadee2014a, Chyong2014} we find that the cost functions are chosen affine, except for production for which a quadratic and a logarithmic term is added to represent the sharp increase in cost close to the capacity limit. (The quadratic and the logarithmic term are both convex.) We adopt these cost structures, except that we omit the logarithmic term in the production cost functions. This allows the problem to be represented as linear complementarity problem while maintaining strict convexity in the production cost function. The derivatives of our cost functions read
\begin{align}
\frac{dc^{Z^\prime}_{zt}(s^{Z^\prime}_{zt})}{ds^{Z^\prime}_{zt}} &= \LINC^{Z^\prime}_{zt},\, Z^\prime \in Z \setminus \{P\} \quad \forall z,t, \\
\frac{dc^P_{nt}(s^P_{nt})}{ds^P_{nt}} &= \LINC^P_{nt} + \QUAC^P_{nt} s^P_{nt} \\
															 &= \LINC^P_{nt} + \QUAC^P_{nt} q^P_{fnt} \quad \forall n,t,
\end{align}
where $\LINC^Z_{zt}$ and $\QUAC^Z_{zt}$ are the linear and quadratic cost function terms for service $Z$ located at node/arc $z$ in period $t$. Note that the last equality holds because every trader has only one associated producer (Section \ref{sec:P1_GeneralSetting}).

We state the full problem including all mentioned changes in \ref{app:P1_ModelEquations}. After dividing the left-hand side of Equation \eqref{eqn:P1_dpiC} by the slope of the inverse demand curve $\SLP^C_{nt}$, the problem can be brought into the following form:
\begin{subequations} \label{eqn:P1_LCP_basic}
\begin{align} 
Mx+b &\succeq 0, \label{eqn:P1_Mxplusq}\\
x &\succeq 0, \label{eqn:P1_xgeq0}\\
x^T \left(Mx+b \right) &= 0, \label{eqn:P1_perp} \\
M:= \begin{bmatrix}
D		&-E^T	&-F^T	&-G^T \\
E		&0		&0		&0 \\
F		&0		&0		&0 \\
G		&0		&0		&H \end{bmatrix},
\quad &x:=\begin{bmatrix} q\\\alpha\\ \phi \\ \lambda \end{bmatrix},
\quad b:=\begin{bmatrix} b_q\\b_\alpha\\ b_\phi \\ b_{\lambda} \end{bmatrix},
\end{align}
\end{subequations}
where $M \in \mathbb{R}^{p \times p}$, $x \in \mathbb{R}^p$, and $b \in \mathbb{R}^p$. As usual, $\succeq, \preceq, \succ, \prec$ are interpreted component-wise. The components of $x$ are defined as
\begin{equation} \label{eqn:P1_x}
q = \begin{bmatrix} q^{P} \\ q^{I} \\ q^{X} \\ q^{A} \\ q^{B} \\ q^{C} \end{bmatrix},\quad
\alpha = \begin{bmatrix} \alpha^Z \\ \alpha^{ZT} \end{bmatrix}, \quad
\phi = \begin{bmatrix} \phi^{N} \\ \phi^{S}\end{bmatrix}, \quad
\lambda = \begin{bmatrix} \lambda^{C}\end{bmatrix},
\end{equation} 
where each sub-component of $q, \alpha, \phi$ and $\lambda$ is again a vector containing all variables of its respective type, for example the vector $q^P$ contains all the produced volumes $q^P_{fnt}$ and is of appropriate dimension. The first block row of $M$ and $b_q$ corresponds to Equations \eqref{eqn:P1_dqP} - \eqref{eqn:P1_dqC} and contains the stationarity conditions of the traders' optimization problems (with substituted service prices $\lambda^Z_{zt}$). $Eq+b_\alpha$ corresponds to Equations \eqref{eqn:P1_dalphaZ} and \eqref{eqn:P1_dbetaZ}, and consists of the primal feasibility constraints of the service providers' problems (with substituted decision variables $s^Z_{zt}$). $Fq+b_\phi$ corresponds to Equations \eqref{eqn:P1_dphiN} and \eqref{eqn:P1_dphiS}, comprising the primal feasibility constraints of the traders' problems. $Gq+H\lambda+b_{\lambda}$ corresponds to the wholesale market clearing conditions contained in Equations \eqref{eqn:P1_dpiC}. Note that the inverse demand functions $\Lambda^{C}_{nt}(\cdot)$ have to be affine for all $n$ and $t$ to represent Model \eqref{eqn:P1_FullProblem} by Equations \eqref{eqn:P1_LCP_basic}. Some of the sub-matrices of $M$ and $b$ have special properties: $D$ is diagonal and positive semi-definite, $E \succeq 0$, $G \preceq 0$, $H$ is diagonal and positive definite, $b_q \succeq 0$, $b_\alpha \succ 0$, $b_\phi = 0$, and $b_\lambda \prec 0$. 

\subsection{Existence of a solution} \label{sec:P1_ExistenceOfASolution}
According to \citet[Theorem 3.1.2]{Cottle1992}, a solution to Problem \eqref{eqn:P1_LCP_basic} exists if
\begin{enumerate}[(i)]
	\item $M$ is positive semi-definite,
	\item a vector $x$ exists satisfying Equations \eqref{eqn:P1_Mxplusq} and \eqref{eqn:P1_xgeq0}.
\end{enumerate}

Since
\begin{align} \label{eqn:P1_M_psd}
x^T M x 
&= \begin{bmatrix} q& \alpha & \phi & \lambda \end{bmatrix} \begin{bmatrix} D&E&F&G \\ -E^T&0&0&0 \\ -F^T&0&0&0 \\ -G^T&0&0&H \end{bmatrix} \begin{bmatrix} q\\ \alpha \\\phi \\\lambda \end{bmatrix} \\
&= q^T D q + \lambda^T H \lambda,
\end{align}
Condition (i) is satisfied if $D$ and $H$ are positive (semi-) definite, which is indeed the case. 

Condition (ii) is proven by construction. We begin with $x_0 = 0$, and by exploiting the structure of $M$ and $b$ we subsequently alter the components of $x$ until we find $x$ satisfying Equations \eqref{eqn:P1_Mxplusq} and \eqref{eqn:P1_xgeq0}.
\paragraph{Iteration 1} $x_0$ is not feasible, because $-G^Tq_0 + H\lambda_0 + b_\lambda = b_\lambda \prec 0$. Therefore, we set $\lambda_1 = -H^{-1}b_\lambda \implies \lambda^C_{nt} = SLP^C_{nt} \cdot \frac{\INT^C_{nt}}{SLP^C_{nt}} = \INT^C_{nt} \geq 0$. Hence, the components of $x_1$ are:
\begin{equation*}
q_1=0,\quad \alpha_1=0, \quad \phi_1=0, \quad \lambda_1 = \INT^C,
\end{equation*}
where $\INT^C$ is the vector containing the intercepts of the inverse demand curves $\INT^C_{nt}$ at the positions corresponding to the indices $n$, $t$ of the wholesale prices $\lambda \equiv \lambda^C$.
\paragraph{Iteration 2} We now consider the first block-row of inequalities: $Dq_1 +E\alpha_1 +F\phi_1 +G\lambda_1 + b_q = G \cdot \INT^C + b_q \not\succeq 0$, because we cannot guarantee $b_q \succeq -G \cdot \INT^C$. We set $\phi^{N}_2 := \INT^C$ (again, $\INT^C$ is a vector containing $\INT^C_{nt}$ at the corresponding indices), which gives us
\begin{equation*}
q_2^T=0,\quad \alpha_2^T=0,\quad \phi_2^N=\INT^C,\quad \phi_2^S=0,\quad \lambda_2 =\INT^C.
\end{equation*}
\paragraph{Iteration 3} Since $F$ has both positive and negative entries, we have $F\phi_2 \preceq 0$ for some lines of the first block-row. Since $E \succeq 0$, and at least one entry on each conflicting line is greater than zero (this can be deduced from the Problem formulation \eqref{eqn:P1_NoLOSSnoTauNoS}), we can compensate for $F\phi_2 \preceq 0$ by setting the corresponding components of $\alpha$ sufficiently large: $\alpha_3 := \INT^C \implies E \alpha \succeq \lvert-F\phi_2 \rvert$. This leads to a feasible $x_3$:
\begin{equation*}
q_3=0, \quad \alpha_3 = \INT^C, \quad \phi^N_3 = \INT^C, \quad \phi^S_3 = 0, \quad \lambda_3 = \INT^C,
\end{equation*}
which completes the proof of Condition (ii).

\subsection{Characterization of the solution space} \label{sec:P1_SolSpace}
As stated in \citet[Theorem 3.3.7]{Cottle1992}, Problem \eqref{eqn:P1_LCP_basic} attains a unique solution for all vectors $b$ if and only if $M$ is a \textbf{\textit{P}}-matrix\footnote{A matrix is a \textbf{\textit{P}}-matrix if it has square form and all its principle minors are positive. This is the most general class of matrices for which the linear complementarity program \eqref{eqn:P1_LCP_basic} attains a unique solution for all vectors $b$. We refer to \citet[Chapter 3.3]{Cottle1992} for an in-depth discussion on \textbf{\textit{P}}-matrices.}. This is clearly not the case as some diagonal elements $M_{ii}=0$ independent of the parameter value choices. However, given $M+M^T$ is diagonal, we can guarantee uniqueness of the solution $x$ in component $i$ if the corresponding element $M_{ii}>0$; see \ref{app:P1_ComponentwiseUniqueness} for a proof. The converse is not necessarily true: $x_i$ may be unique even if $M_{ii}=0$ because the component is sufficiently constrained by $M$ and $b$.

Consequently, results are unique in at least the produced volume $q^P$, the sold volume $q^{C}_{fnt}$ if the corresponding market power parameter $\theta^C_{fnt}>0$, and the wholesale price $\lambda^C$. Furthermore, from Equation \eqref{eqn:P1_MarketClearingEndConsumer} and the properties of the inverse demand function \eqref{eqn:P1_LambdaC}-\eqref{eqn:P1_eta_nt}, we can deduce that the total consumption in a node and time period $\sum \limits_{f \in \mathcal{F}(n)} q^{C}_{fnt}$ is unique. By combining these results, we can conclude that the total sales of perfectly competitive traders in a node and time period $\sum \limits_{f \in \{\mathcal{F}(n) | \theta_{fnt} = 0\}} q^{C}_{fnt}$ are unique, including the special case of $q^{C}_{fnt}$ being unique independently of $\theta_{fnt}$ if there is only one trader $f$ in the market in node $n$ and time period $t$. Finally, unique produced volumes $q^P$ imply that the sold volume $q^{C}_{fnt}$ is unique if trader $f$ only sells gas in one node and one time period. 

To determine uniqueness of $x_i$ for all $i$ for specific instances of $M$ and $b$ we analyze the solution space of the problem: Since $M$ is positive semi-definite, the solution space $\mathbb{S}$ is a convex polyhedron \citep[Theorems 3.1.7 and 3.1.8]{Cottle1992} of the form
\begin{equation} \label{eqn:P1_setS}
\mathbb{S} = \{ x \in \mathbb{R}^p_{+} : Mx+b \succeq 0, b^T(x-\hat{x}) = 0, (M+M^T)(x-\hat{x}) = 0 \},
\end{equation}
where $\hat{x}$ is an arbitrary solution. We have shown in Section \ref{sec:P1_ExistenceOfASolution} that a solution $\hat{x}$ always exists, hence the solution space can always be described by $\mathbb{S}$.

We are now interested in the extreme points of $\mathbb{S}$, particularly those attaining maximum or minimum values in (at least) one component: Let $x^{i,+}, x^{i,-} \in \mathbb{S}$ be vectors (extreme points) attaining the maximum and minimum values of $\mathbb{S}$ in component $i$. 
Clearly, the solution to Problem \eqref{eqn:P1_LCP_basic} is unique in component $i$ if and only if $x^{i,+}=x^{i,-}$.

To compute $x^{i,+}$ and $x^{i,-}$ for all $i$ for a specific $M$ and $b$, we first calculate an arbitrary solution $\hat{x}$ by solving Problem \eqref{eqn:P1_LCP_basic}. Second, we solve $2p$ Linear Programs of the form 
\begin{equation} \label{eqn:P1_LP}
\min \limits_{x \in \mathbb{S}} \left(c^{i,\pm}\right)^T x,
\end{equation}
where $c^{i,\pm} \in \mathbb{R}^p$ is $\pm 1$ in component $i$ and 0 in all other components. 

As noted in Section \ref{sec:P1_Introduction}, lack of uniqueness may have profound implications for the interpretation of the solutions to Problem \eqref{eqn:P1_LCP_basic}. These implications will be discussed in detail in Section \ref{sec:P1_Discussion}, after we provide an example.

\section{Numerical example} \label{sec:P1_NumericalExample}
In this section, we demonstrate the multiplicity of solutions via a numerical example. Based on the data sources listed in \ref{app:P1_AvailableDataSources}, we set up Model \eqref{eqn:P1_NoLOSSnoTauNoS} to replicate the global gas market in 2012. The model comprises $50$ nodes and $291$ arcs, distinguishes the summer and winter half-year, and covers \mbox{97\%} of the global gas production and consumption. We refrain to go into the details of the modeling and calibration process, as any parameters fulfilling the criteria listed in the previous sections would qualitatively lead to the same outcome and therefore the exact values of the parameters are not relevant for the demonstration of the solution properties. The reason for underlying the simulations with real data is merely to generate results which are in the order of magnitude typically seen in gas market analyses. 

We simulated two cases, a \ac{BC} and a \ac{CF}, allowing for the comparison of (i) multiple solutions in one case, and (ii) the solution spaces of the two cases. In \ac{BC}, we increased the market power parameter $\theta^C_{fnt}$ to $0.01$ if the calibrated value was less than $0.01$, while leaving all $\theta^C_{fnt} \in [0.01,1]$ at their calibrated values. This allowed testing whether $q^C_{fnt}$ attains a unique result if $\theta^C_{fnt}>0$ (Section \ref{sec:P1_SolSpace}), while avoiding conflicts with the solver's tolerances\footnote{We used the standard settings of the Gurobi solver for Quadratic Programs (Version 5.6.3).} on one hand, and keeping the impact on the underlying market behavior low on the other hand. In \ac{CF}, we changed the traders' market power parameters $\theta^C_{fnt}=0$ for all $f,n,t$. \ac{CF} corresponds to the (hypothetical) situation of perfect competition among all traders in all markets and time periods which is often taken as benchmark in case studies. We followed the procedure described in Section \ref{sec:P1_SolSpace} to determine the solution spaces of \ac{BC} and \ac{CF}. The results are shown in Tables \ref{tab:P1_CaseUniqueness} and \ref{tab:P1_CaseSolSpace}.

Table \ref{tab:P1_CaseUniqueness} provides the maximum differences between all possible market equilibria in \ac{BC} and \ac{CF}, grouped by similar variables. We observe that the solution is unique in the produced volumes $q^P$, the dual variables associated with the capacity constraints of the services $\alpha$, the wholesale prices $\lambda^C$, the volume flows contracted by the service providers $s^Z$, and the service prices $\lambda^Z$. However, the solution is not unique in the other volume flows (corresponding to the remaining components of $q$), and the dual variables associated with the volume balances of the traders $\phi$. Compared to \ac{BC}, the maximum differences in \ac{CF} for the \ac{LNG} volumes $q^B$ and sold volumes $q^C$ switch from zero to non-zero, while all other components remain unique/non-unique. Also, in \ac{CF} the maximum differences of the volumes $q$ are larger. The maximum difference in $\phi$ is also shown for the sake of completeness, although this result is not of interest in interpretation.

\begin{table}[!htb] 
\centering
\caption{Maximum difference between all possible market equilibria in \ac{BC} and \ac{CF}, grouped by variables: $\max \limits_{i \in J} \left(x^{i,+}_i-x^{i,-}_i\right)$, where $J$ is the index set of the group of variables. All volume flows $q$ and $s^Z$ are given in million cubic meters per day (\unitfrac{mcm}{d}). All prices $\lambda^C$ and $\lambda^Z$, as well as the dual variables $\alpha$ and $\phi$, are given in thousand Euros per million cubic meters (\unitfrac{k\euro}{mcm}). The table also indicates for which variable the maximum difference is attained (Attained by trader $f$, node $n$/arc $nm$, period $t$), and gives the maximum value attained per group of variables over all simulations (Max. value attained BC/CF). RU: Russia. UA: Ukraine. AZ: Azerbaijan. NAm: North America. CeA: Central Asia, including Kazakhstan, Turkmenistan and Uzbekistan. JP: Japan. DE: Germany. ap: April-September. oc: October-March.}
\label{tab:P1_CaseUniqueness}
\begin{tabular}{c|r:p{0.3cm}p{0.3cm}p{0.3cm}p{0.3cm}|r:p{0.3cm}p{0.3cm}p{0.3cm}p{0.3cm}|r}
\hline 
Group 			& \multicolumn{1}{c:}{Max.} 				& \multicolumn{4}{c|}{Attained by} 																																			& \multicolumn{1}{c:}{Max.} 			& \multicolumn{4}{c|}{Attained by} 																																		& \multicolumn{1}{c}{Max.}\\
 of					& \multicolumn{1}{c:}{diff.}				& & & &																																																	& \multicolumn{1}{c:}{diff.} 			& & & &																																																& \multicolumn{1}{c}{value}  \\ 
var.				& \multicolumn{1}{c:}{in} 					& & & &																																																	& \multicolumn{1}{c:}{in} 				& & &	&																																																& \multicolumn{1}{c}{attained} \\ 
						&	\multicolumn{1}{c:}{\ac{BC}}			& \multicolumn{1}{c}{$f$} &\multicolumn{1}{c}{$n$} & \multicolumn{1}{c}{$m$} &\multicolumn{1}{c|}{$t$}	& \multicolumn{1}{c:}{\ac{CF}}		& \multicolumn{1}{c}{$f$} &\multicolumn{1}{c}{$n$} & \multicolumn{1}{c}{$m$} &\multicolumn{1}{c|}{$t$}& \multicolumn{1}{c}{in \ac{BC}/\ac{CF}}\\ \hline
$q^{P}$ 		& \textbf{0}  & & & &															& \textbf{0} 	& &	& &				 											& 2530\\ 
$q^{I}$			& \textbf{21} & RU & RU &  & ap 									& \textbf{137}& RU & UA &  & ap										& 335\\ 
$q^{X}$ 		& \textbf{21} & RU & RU & - & oc 									& \textbf{135}& RU & UA & - & oc									& 332\\ 
$q^{A}$ 		& \textbf{21} & RU & RU & UA & ap									& \textbf{171}& CeA & RU & UA & oc								& 397\\ 
$q^{B}$ 		& \textbf{0}  & & & &  														& \textbf{43} & RU & RU & JP & ap									& 267\\ 
$q^{C}$ 		& \textbf{0}  & & & &  														& \textbf{186}& RU & DE &  & oc										& 2656\\ 
$\alpha$ 		& \textbf{0}	& & & &  														& \textbf{0}	& & & & 														& 237\\ 
$\phi$ 			& \textbf{488}& AZ & NAm &  & 										& \textbf{340}& AZ & NAm & &											& 653\\ 
$\lambda^C$	& \textbf{0} 					& & & &  														&	\textbf{0}	 					& & & &  														& 490\\ 
$s^Z$ 			& \textbf{0} 					& & & &  														&	\textbf{0}	 					& & & & 														& 2530\\ 
$\lambda^Z$	& \textbf{0} 					& & & &  														& \textbf{0}	 					& & & & 														& 419\\ \hline 
\end{tabular}
\end{table}

Table \ref{tab:P1_CaseSolSpace} compares the minimum and maximum values achieved, $x^{i,-}$ and $x^{i,+}$, in \ac{BC} and \ac{CF} for selected components $i$. The examples show that the market equilibrium in \ac{CF} can be higher or lower than in \ac{BC} for the sold volumes $q^C$. For the injected, extracted and piped volumes $q^I$, $q^X$, and $q^A$, which are ambiguous in \ac{BC} and \ac{CF}, the examples demonstrate that the solution intervals overlap. Thus, in none of the presented examples a prediction of the market equilibrium is possible when moving from \ac{BC} to \ac{CF}.

\begin{table}[!htb]
\centering
\caption{Comparison of components $i$ of exemplary extreme points $x^{i,+}$, $x^{i,-}$ of the solution spaces of \ac{BC} and \ac{CF}. $\underline{BC}$: Minimum achieved value in component $i$ in \ac{BC}. $\overline{BC}$: Maximum achieved value in component $i$ in \ac{BC}. $\underline{CF}$: Minimum achieved value in component $i$ in \ac{CF}. $\overline{CF}$: Maximum achieved value in component $i$ in \ac{CF}. $f$: trader. $n$: (starting) node in which transaction takes place. $m$: ending node in which transaction takes place. $t$: time period. NO: Norway. GB: Great Britain. RU: Russia. CH: Switzerland. IT: Italy. NL: the Netherlands. DE: Germany. ap: April-September. oc: October-March. All figures are given in million cubic meters per day (\unitfrac{mcm}{d}).}
 \label{tab:P1_CaseSolSpace}
\begin{tabular}{l|cccc|rrrr}
\hline
Variable & $f$ & $n$& $m$& $t$& $\underline{BC}$ & $\overline{BC}$& $\underline{CF}$ & $\overline{CF}$\\ \hline
$q^{C}_{fnt}$ 	& NO & GB & - & ap & 37.6 & 37.6 & 24.5 & 45.5 \\ 
$q^{C}_{fnt}$ 	& RU & CH & - & ap & 0.77 & 0.77 & 0 		& 5.06 \\
$q^{I}_{fnt}$ 	& NO & IT & - & ap & 7.01 & 18.9 & 0		& 9.77 \\
$q^{X}_{fnt}$ 	& NL & IT & -	& oc & 4.11 & 16.5 & 0 		& 9.68 \\
$q^{A}_{fnmt}$ 	& NO & DE & CH& oc & 12.1 & 24.3 & 0 		& 14.6 \\
\hline
\end{tabular}
\end{table}

\section{Discussion} \label{sec:P1_Discussion}
The theoretical results derived in Section \ref{sec:P1_ExistanceAndUniqueness} and the numerical example presented in Section \ref{sec:P1_NumericalExample} provide deep insight in the mathematical structure of Problem \eqref{eqn:P1_LCP_basic} and allow drawing conclusions about the uniqueness of the individual components of the solution. In particular, the numerical example demonstrates how large the difference between two possible solutions can be in a practical case and thereby emphasizes the importance of analyzing the entire solution space when carrying out a case study. In the following, we discuss various aspects of our findings in detail.

As Table \ref{tab:P1_CaseUniqueness} indicates, the solution is unique in the produced volumes $q^P$, which is in line with our theoretical findings in Section \ref{sec:P1_SolSpace} and \ref{app:P1_ComponentwiseUniqueness}. Furthermore, we proved in Section \ref{sec:P1_SolSpace} that the sold volume $q^C_{fnt}$ is unique if the corresponding market power parameter $\theta^C_{fnt}>0$; this result was also confirmed by the numerical example (Table \ref{tab:P1_CaseUniqueness}). From \ac{CF} in Table \ref{tab:P1_CaseUniqueness}, we see that the ambiguity can be very large if the market power parameter $\theta^C_{fnt} = 0$: Some sold volumes $q^C_{fnt}$ attain solutions as different as \unitfrac[186]{mcm}{d}. As a comparison, Germany, the largest European gas consumer, consumed \unitfrac[293]{mcm}{d} during the winter season 2012. Hence, the sold volumes could not be predicted in general for traders behaving perfectly competitive.

Fortunately, this is not required for decision-making in most cases; The key variable in a situation with perfectly competitive traders is the wholesale market price $\lambda^C_{nt}$, which is unique independently of the market power parameters $\theta^C_{fnt}$ of the traders in node $n$ and time period $t$ (Table \ref{tab:P1_CaseUniqueness}). However, if some perfectly competitive traders have high market shares $\frac{q^C_{fnt}}{\sum \limits_{f \in \mathcal{F}(n)} q^{C}_{fnt}}$, countries might still be concerned: Traders might develop strategic behavior over time, and such a transition may induce large changes in the volume flows, as illustrated by our example (Table \ref{tab:P1_CaseSolSpace}).

Alternatively, by choosing the market power parameter $\theta^C_{fnt}>0$ for all traders $f$ in all nodes $n$ and time periods $t$ we can guarantee uniqueness of the volumes sold $q^C$. However, the underlying model assumptions leading to this specific, unique $q^C$ should be studied carefully: first, $\theta^C_{fnt}>0$ for all $f$, $n$, $t$, implies that all traders in the model exert market power in all markets and time periods to some extent. We argue that this is plausible for a real world case, because real markets are imperfect and therefore some market power can always be exerted by traders. Second, a market power parameter $\theta^C_{fnt}>0$ drives trader $f$ to shift some of its sold volumes away from market $n$ in period $t$. This originates from the fact that $\theta^C_{fnt}$ is proportional to a cost on $\left(q^C_{fnt}\right)^2$. This becomes obvious when the inverse demand function $\Lambda^{C}_{nt}(\cdot)$ is plugged into the Revenue function \eqref{eqn:P1_TraderRevenue} of the trader. If $\theta^C_{fnt}>0$ for all $f$, $n$, $t$, all traders diversify their consumers. Since we indeed observe that traders diversify their consumers in the real world, for instance to reduce risks associated with the counter-party, we conclude that $\theta^C_{fnt}>0$ for all $f$, $n$, $t$ can be justified from an economic perspective and the resulting volume sold $q^C_{fnt}$ is well defensible. Hence, in practice it might be preferable to assign a small positive value to the market power parameter $\theta^C_{fnt}$ of perfectly competitive traders, since the resulting volumes sold $q^C_{fnt}$ are unambiguous and plausible, and the exploration of the solution space can be circumvented.

The remaining volume flows $q^{Z^{\prime \prime}}$, $Z^{\prime \prime} \in \{I,X,A,B\}$ do not attain unique solutions in the simulations (Table \ref{tab:P1_CaseUniqueness}). On one hand, we could again achieve uniqueness by enforcing the elements in the diagonal of $D$, $D_{ii}$, to be larger than zero for those $i$ for which $x_i = q^{Z^{\prime \prime}}_{fzt}$. This would simplify the interpretation of the results but is not justifiable from an economic perspective: $D_{ii} > 0$ implies that the corresponding trader $f$ diversifies service $Z^{\prime \prime}$ in $z$ and $t$, because $D_{ii}$ is a cost on $\left(q^{Z^{\prime \prime}}_{fzt}\right)^2$. This is slightly more difficult to see from the model equations \eqref{eqn:P1_NoLOSSnoTauNoS}, because no such term on $\left(q^{Z^{\prime \prime}}_{fzt}\right)^2$ existed in the original problem \eqref{eqn:P1_FullProblem}. However, by adding $D_{ii} q^{Z^{\prime \prime}}_{fzt}$ to the corresponding inequality in Model \eqref{eqn:P1_NoLOSSnoTauNoS}, the analogy to the term $-\theta^C_{fnt} \SLP^C_{nt} q^C_{fnt}$ in Equation \eqref{eqn:P1_dqC} becomes obvious, which in turn was derived from the Revenue function \eqref{eqn:P1_TraderRevenue} and can be interpreted as a quadratic cost term as discussed above. However, traders in the real world often receive quantity discounts from service providers and therefore rarely diversify. This would correspond to $D_{ii}<0$, which would lead to an overall non-convex problem and is therefore not considered here.

On the other hand, in practice multiple market equilibria can indeed be in line with a given situation, as illustrated by the following example: assume multiple traders which are moving gas from one node to another. The two nodes are interconnected by two different paths, of which the path with lower marginal costs is congested. In equilibrium, the costs for moving gas from one node to the other one are equal for all traders regardless of the path they are using, because the congestion fee charged on the cheaper path is equal to the difference in marginal costs of the two paths. Hence, which trader will eventually use which path in reality is due to factors that cannot at first glance be related to the model, such as the trader's preferences, the point in time the capacity is booked, historical reasons, etc. We therefore recommend not to alter the model structure and instead analyze the solution space as a whole when the interpretation of these variables is of interest.

For all other variables we cannot guarantee uniqueness although they are unique in the presented example, as Table \ref{tab:P1_CaseUniqueness} indicates. This is particularly disadvantageous when it comes to the volumes processed by the service providers $s^Z$, and the corresponding service prices $\lambda^Z$, since these variables are valuable for planning infrastructure expansions. In fact, we can easily find counterexamples showing that $s^Z$ and $\lambda^Z$ are not unique in general: Assume $2$ uncongested paths with identical start and end nodes running through different intermediate nodes. If the marginal costs for gas transportation are equal on both paths, then multiple optimal $s^A$ exist. Now assume a path through $3$ nodes, where the second node is only connected to the start and end node. If both arcs are congested, $\lambda^A$ is not unique; only the sum of all $\lambda^A$ of these two arcs is unique \citep[see also][p.~231]{Egging2010a}. Similar examples can also be found for the other services. This said, when models reproduce a real world network, situations with non-unique $s^Z$ and $\lambda^Z$ are rare, and therefore the interpretation of these variables is often facilitated. Nonetheless, the entire solution space has to be explored to guarantee uniqueness.

Finally, we discuss extensions of our results to more involved models. A first group comprises the stochastic versions of the \ac{WGM} introduced by \citet{Zhuang2008} and \citet{Gabriel2009}. In these models the intercept $\INT^C$ and slope $\SLP^C$ of the inverse demand functions are randomized (within the usual bounds, $\INT^C>0$ and $\SLP^C<0$), and the possible realizations are reflected by a scenario tree. Since our results are valid for any $\INT^C>0$ and $\SLP^C<0$, we can conclude that the decision variables of the individual tree nodes suffer from the ambiguity described above, and therefore, the overall outcome of the model is ambiguous as well.

Another popular extension is the introduction of convex logarithmic terms in the production cost functions. Since the uniqueness-result derived in \ref{app:P1_ComponentwiseUniqueness} holds for any strictly convex production cost function, we conclude that the produced volume $q^P$ is also unique for this group of models. Furthermore, the solution space $\mathbb{S}$ is polyhedral and can be described similarly as before (Section \ref{sec:P1_SolSpace}): First, we solve the non-linear complementarity problem $0\preceq K(x)\perp x \succeq 0$ and obtain a solution $\hat{x}$. Second, we express $\mathbb{S}$ by means of $K(x)$, $x$ and $\hat{x}$: We set $q^P = \hat{q}^P$ and insert the numerical values into $K(x)$, which leads to an affine $K(x)=\bar{M}x+\bar{b}$. Finally, we solve the Linear Problem \eqref{eqn:P1_LP} and obtain the extreme points of interest $x^{i,+}$, $x^{i,-} \in \mathbb{S}$.

A further example are the introduced and then dropped constraints \eqref{eqn:P1_UpperCLimits}-\eqref{eqn:P1_LowerZLimits} on the volume flows $q^C$ and $q^Z$. They would be reflected by additional terms in the KKT conditions \eqref{eqn:P1_FullProblem}. However, the fundamental properties of matrix $M$, particularly positive semi-definiteness and diagonality of $M+M^T$,  remain unchanged, and thus the derived results regarding existence and uniqueness of the solutions extend to models including constraints \eqref{eqn:P1_UpperCLimits}-\eqref{eqn:P1_LowerZLimits}.
\section{Conclusions} \label{sec:P1_Conclusions}
In this paper we introduced a SPE-CV model for the natural gas market. We specified affine costs for the service providers, quadratic costs for gas production, and affine inverse demand, and hence represented our model as linear complementarity problem. We proved that a solution to the given problem always exists and derived criteria for uniqueness of the individual components. Finally, we presented a numerical example confirming our theoretical findings and illustrating the main difficulties for the interpretation of the results. 

Our findings indicate that most of the relevant variables for decision-making are unique in the presented model setting. However, one has to be particularly careful when drawing conclusions about the volume flows between the traders contract with the service providers, and the traders sell to the consumers if the traders do not exert market power, since in these cases the model allows multiple market equilibria for a fixed set of parameter values. The obtained set of solutions $\mathbb{S}$ isolates the possible market equilibria, but which of the equilibria will materialize in practice cannot be predicted without further specifying the traders' preferences.

The situation is even more difficult when one tries to compare market equilibria computed under different sets of parameter values: different results for the above quantities may not be due exclusively to the different parameter values, but also to non-uniqueness of the solution for each set of parameter values. In the worst case, the direction of change of some variables is ambiguous, and therefore no trend can be predicted from the simulation results. 

Fortunately, we can guarantee uniqueness in the gas sales of individual traders to consumers $q^{C}$, perhaps the most important non-unique variable in the model, by assuming that all traders exert at least some market power. Moreover, we showed that a trader exerting market power over a group of consumers is mathematically equivalent to a trader diversifying its consumers. Since most of the traders in the real world aim for some consumer diversification, we conclude that this is a practicable approach to enforce uniqueness in $q^{C}$, and that the computed $q^{C}$ is sensible from an economic perspective. 

Since our findings are linked to the mathematical structure of the model, they are not limited to gas markets in particular. Instead, they apply to any other oligopolistic market which is representable by the introduced model formulation. In fact, the results also hold if the presented structure is part of a larger model, for instance one including multiple energy carriers as proposed by \citet{Huppmann2014}, and therefore may be useful in a large field of applications.

\section*{Acknowledgements}
\addcontentsline{toc}{section}{Acknowledgements}
We would like to thank our colleague Dr.~Peyman Esfahani for the helpful discussions on uniqueness of solutions of linear complementarity programs, and Prof.~Ruud Egging and the four anonymous reviewers for their valuable comments.

\clearpage
\appendix
\section{Notation} \label{app:Notation}
\begin{table}[htbp]
\centering
\caption{This table introduces the nomenclature concerning service providers, traders and consumers.}
 \label{tab:P1_ServiceProviders}
\def\arraystretch{1.2}
\begin{tabular}{p{0.15\textwidth}p{0.75\textwidth}}
\hline 
\multicolumn{2}{l}{\textbf{Service providers, traders and consumers}} \\
\hline 
$A_{nm}$& Transmission system operator of pipeline $nm$\\
$B_{nm}$& Shipping company transporting \ac{LNG} from $n$ to $m$\\
$C_n$& Consumer at node $n$\\
$I_n$& Storage operator injecting gas at node $n$\\
$L_n$& Liquefaction plant operator at node $n$\\
$P_n$& Gas producing company at node $n$\\
$R_n$& regasification plant operator at node $n$\\
$S_n$& Storage operator at node $n$\\
$F_n$& The trader associated with producer at \mbox{node $n$}\\
$X_n$& Storage operator extracting gas at node $n$\\
$Z_z$& Placeholder for a service provider ($P_n$, $I_n$, $X_n$, $L_n$, $R_n$, $A_{nm}$, $B_{nm}$) at node $n$ / arc $nm$\\
\hline
\end{tabular}
\end{table}

\begin{longtable}[htbp]{p{0.33\textwidth}p{0.57\textwidth}}
\caption{This table introduces all sets used for the mathematical description of the model.}
\label{tab:P1_Sets} \\
\hline 
\multicolumn{2}{l}{\textbf{Sets}}  \\
\hline 
\endfirsthead
\multicolumn{2}{c}{\begin{footnotesize}\tablename\ \thetable\ -- \textit{Continued from previous page}\end{footnotesize}} \\
\hline
\multicolumn{2}{l}{\textbf{Sets}}  \\
\hline
\endhead
\hline \multicolumn{2}{r}{\begin{footnotesize}\textit{Continued on next page}\end{footnotesize}} \\
\endfoot
\hline
\endlastfoot
$t \in \mathcal{T} = \{T_1, \ldots, T_{\bar{t}}\}$& A time period $t$ in the set $\mathcal{T}$ of all periods of a year\\
$n, m \in \mathcal{N} = \{N_1, \ldots, N_{\bar{n}}\}$& Nodes $n,m$ in the set $\mathcal{N}$ of all nodes \\ 
$f \in \mathcal{F} = \{F_1, \ldots, F_{\bar{n}}\}$& A trader $f$ in the set $\mathcal{F}$ of all traders\\
$z \in \mathcal{Z}$ & A node/arc element from the \mbox{set $\mathcal{Z}$} \\ \hline
$\mathcal{A} \subset \mathcal{N} \times \mathcal{N}$& Set of arcs connecting $2$ nodes by pipeline					\\
$\mathcal{B} \subset \mathcal{N} \times \mathcal{N}$ & Set of arcs connecting $2$ nodes by ship					\\	
$\mathcal{C} \subseteq \mathcal{N}$ & Set of nodes at which a consumer is active \\
$\mathcal{I} \subseteq \mathcal{N}$ & Set of nodes at which storage injection	is possible					\\
$\mathcal{L} \subseteq \mathcal{N}$ & Set of nodes at which a liquefaction terminal operator is active				\\	
$\mathcal{P} \subseteq \mathcal{N}$ & Set of nodes at which a gas producer is active	\\							
$\mathcal{R} \subseteq \mathcal{N}$ & Set of nodes at which a regasification terminal operator is active		\\
$\mathcal{X} \subseteq \mathcal{N}$ & Set of nodes at which storage extraction is possible						\\
$\mathcal{Z} \in \{\mathcal{P},\mathcal{L},\mathcal{B},\mathcal{R},\mathcal{A},\mathcal{I},\mathcal{X}\}$ & Placeholder for the set of nodes/arcs at which a type of service provider is active \\ \hline
$\mathcal{A}(n) \subseteq \mathcal{N} \setminus \{n\}$ & Set of nodes which are connected to $n$ by pipeline \\
$\mathcal{B}(n) \subseteq \mathcal{N} \setminus \{n\}$ & Set of nodes which are connected to $n$ by ship \\
$\mathcal{C}(f) \subseteq \mathcal{N}$& The set of all nodes with consumers which are reachable by \mbox{trader $f$}\\
$\mathcal{N}(f) \subseteq \mathcal{N}$& The set of all nodes which are reachable by \mbox{trader $f$}\\
$\mathcal{F}(z)$& The set of all traders active at node/arc $z$\\
$\mathcal{Z}(f)$& The set of all nodes/arcs in which service $Z$ is active and are reachable by \mbox{trader $f$}\\ \hline
$\mathbb{S}$ & The solution space to Problem \eqref{eqn:P1_LCP_basic} \\
\hline
\end{longtable}

\begin{table}[htbp]
\centering
\caption{The \textit{parameters} are generally described by capital Roman letters, and occasionally by lower-case Greek letters to follow conventions. The superscripts indicate whether the parameter is related to a service provider of type $Z \in \{P,L,B,R,A,I,X\}$, a consumer $C$, and whether the all periods are covered ($T$). Subscripts indicate the trader $f$, node/arc $z$, and the period of the year $t$ the parameter is related to.}
 \label{tab:P1_Parameters}
\def\arraystretch{1.2}
\begin{tabular}{p{0.12\textwidth}p{0.78\textwidth}}
\hline 
\multicolumn{2}{l}{\textbf{Parameters}} \\
\hline 
$\overline{\CAP}^Z_{nt}$& Maximum capacity of service $Z$ located at $z$ in \mbox{period $t$}\\
$\overline{\CAP}^{ZT}_{n}$& Maximum capacity of service $Z$ located at $z$ over all \mbox{periods $\mathcal{T}$}\\
$\DMD^C_{nt}$& Reference demand used for construction of the demand curve\\
$\INT^C_{nt}$ & Maximum willingness to pay of consumers at node $n$ in \mbox{period $t$} \\
$\LINC^Z_{zt}$ & Linear cost function term for service $Z$ located at $z$ in \mbox{period	$t$}\\
$\LOSS^Z_z $& Loss factor when using service $Z$  located \mbox{at $z$} \\
$\QUAC^Z_{zt}$ & Quadratic cost function term for service $Z$ located at $z$ in \mbox{period $t$}\\
$\SHR^C_{nt,(\cdot)}$& Reference share in demand of domestic, industry and electricity producing sectors \\
$\SLP^C_{nt}$& Slope of the inverse demand curve of the consumers at node $n$ in \mbox{period $t$}, is assumed strictly negative\\
$\WTP^C_{nt}$& Reference willingness to pay used for construction of the demand curve\\
$\theta^C_{fnt}$& Market power parameter of trader $f$ at node $n$ in \mbox{period $t$}\\
$\eta^C_{nt,(\cdot)}$& Price elasticity of the domestic, industry and electricity producing sectors\\
\hline
\end{tabular}
\end{table}

\begin{table}[htbp]
\centering
\caption{The \textit{variables} are described by lowercase letters. Primal variables are Roman, while dual variables are Greek letters. The superscripts indicate whether the variable is related to a service provider of type $Z \in \{P,L,B,R,A,I,X\}$, a consumer $C$, or a node $N$. Subscripts indicate the trader $f$ the variable corresponds to, at which node/arc $z$ the transaction or service is located, and in which period of the year $t$ it takes place.}
 \label{tab:P1_Variables}
\def\arraystretch{1.2}
\begin{tabular}{p{0.1\textwidth}p{0.8\textwidth}}
\hline 
\multicolumn{2}{l}{\textbf{Variables}} \\
\hline 
$q^{Z}_{fzt}$ & Flow between trader $f$ and service provider $Z$ at node/arc $z$ in \mbox{period $t$} \\
$q^{C}_{fnt}$ & Flow of trader $f$ to consumer $C$ at node $n$ in \mbox{period $t$}\\
$s^{Z}_{zt}$ & Volume flow contracted by service provider $Z$ at node/arc $z$ in \mbox{period $t$} \\
$\alpha^Z_{nt}$ & Congestion fee of service $Z$ at node $n$ in \mbox{period $t$}\\
$\alpha^{ZT}_{n}$ & Congestion fee on annual usage of service $Z$ at \mbox{node $n$} \\
$\phi^N_{fnt}$ & Dual variable of the volume balance of trader $f$ at node $n$ and \mbox{period $t$}\\
$\phi^S_{fn}$ & Dual variable of the annual volume balance of trader $f$ in storage $S$ at \mbox{node $n$}\\
$\lambda^C_{nt}$ & End consumer price at node $n$ in \mbox{period $t$}\\
$\lambda^Z_{nt}$ & Price for utilizing service $Z$ at node $n$ in \mbox{period $t$}\\
\hline
\end{tabular}
\end{table}

\begin{table}[htbp]
\centering
\caption{This table introduces the \textit{functions}. The superscripts indicate whether the function is related to a service provider of type $Z \in \{P,L,B,R,A,I,X\}$, or a consumer $C$. Subscripts indicate at which node $n$ the service/consumer is located, and in which period of the year $t$ the function is valid.}
 \label{tab:P1_Functions}
\def\arraystretch{1.2}
\begin{tabular}{p{0.1\textwidth}p{0.8\textwidth}}
\hline 
\multicolumn{2}{l}{\textbf{Functions}} \\
\hline 
$c^{Z}_{zt}(s^Z_{zt})$ & Cost function of service $Z$ at node/arc $z$ in period $t$.\\
$\Lambda^{C}_{nt}(s^C_{nt})$ & Inverse demand function of consumer $C$ at node $n$ in period $t$.\\
\hline
\end{tabular}
\end{table}

\section{Model equations} \label{app:P1_ModelEquations}
\begin{subequations} \label{eqn:P1_NoLOSSnoTauNoS}
\allowdisplaybreaks
\begin{alignat}{4}
0 &\leq& \, \LINC^P_{nt} + \QUAC^P_{nt} q^P_{fnt} + \alpha^P_{nt} + \alpha^{PT}_{n} - \phi^N_{fnt} & \perp q^P_{fnt} &&\,\geq 0 \quad \forall f,n,t \label{eqn:P1_dqP}\\
0 &\leq&\, \LINC^I_{nt} + \alpha^I_{nt} + \alpha^{IT}_{n} + \phi^N_{fnt} - \phi^S_{fn}& \perp q^I_{fnt} &&\,\geq 0 \quad \forall f,n,t \label{eqn:P1_dqI}\\
0 &\leq&\, \LINC^X_{zt} + \alpha^X_{nt} + \alpha^{XT}_{n} - \phi^N_{fnt} + \phi^S_{fn} & \perp q^X_{fnt} &&\,\geq 0 \quad \forall f,n,t \label{eqn:P1_dqX}\\
0 &\leq&\, \LINC^A_{nmt} + \alpha^A_{nmt} + \alpha^{AT}_{nm} - \phi^N_{fmt} + \phi^N_{fnt} & \perp q^A_{fnmt} &&\,\geq 0 \quad \forall f,n,m,t \label{eqn:P1_dqA}\\
0 &\leq&\, \LINC^L_{nt} + \alpha^L_{nt} + \alpha^{LT}_{n} + \LINC^B_{nmt} + \alpha^B_{nmt} + \alpha^{BT}_{nm} &&& \notag\\
&&\, + \LINC^R_{mt} + \alpha^R_{mt} + \alpha^{RT}_{m} - \phi^N_{fmt} +  \phi^N_{fnt} & \perp q^B_{fnmt} &&\,\geq 0 \quad \forall f,n,m,t \label{eqn:P1_dqB}\\
0 &\leq&\, -\lambda^C_{nt} -\theta^C_{fnt} \SLP^C_{nt} q^C_{fnt} + \phi^N_{fnt}& \perp q^C_{fnt} &&\,\geq 0  \quad \forall f,n,t \label{eqn:P1_dqC}\\
0 &\leq&\, q^{P}_{fnt} + q^{X}_{fnt} + \sum \limits_{m \in \mathcal{A}(n)} q_{fmnt}^{A} + \sum \limits_{m \in \mathcal{B}(n)} q_{fmnt}^{B}&&& \notag \\
&&\,- q^{I}_{fnt} -  q^{C}_{fnt} - \sum \limits_{m \in \mathcal{A}(n)} q_{fnmt}^{A} - \sum \limits_{m \in \mathcal{B}(n)} q_{fnmt}^{B} & \perp \phi^N_{fnt} &&\,\geq 0  \quad \forall f,n,t \label{eqn:P1_dphiN}\\
0 &\leq&\, \sum \limits_{t \in \mathcal{T}} q^{I}_{fnt} - \sum \limits_{t \in \mathcal{T}} q^{X}_{fnt} &\perp \phi^S_{fn} &&\,\geq 0  \quad \forall f,n\label{eqn:P1_dphiS}\\
0 &\leq&\, \overline{\CAP}^Z_{zt} - \sum \limits_{f \in \mathcal{F}(z)} q^{Z}_{fzt} & \perp \alpha^Z_{zt} &&\,\geq 0  \quad \forall z,t \label{eqn:P1_dalphaZ}\\
0 &\leq&\, \overline{\CAP}^{ZD}_{z} - \sum \limits_{t \in \mathcal{T}} \sum \limits_{f \in \mathcal{F}(z)} q^{Z}_{fzt} & \perp \alpha^{ZT}_z &&\,\geq 0  \quad \forall z \label{eqn:P1_dbetaZ}\\
0 &\leq&\, \lambda^C_{nt} - \left( \INT^C_{nt} + \SLP^C_{nt} \sum \limits_{f \in \mathcal{F}(n)} q^{C}_{fnt} \right) &\perp \lambda^C_{nt} &&\,\geq 0  \quad \forall n,t\label{eqn:P1_dpiC}
\end{alignat}
\end{subequations}

\section{Component-wise uniqueness} \label{app:P1_ComponentwiseUniqueness}
Problem \eqref{eqn:P1_LCP_basic}, provided feasible, is equivalent to Problem \eqref{eqn:P1_QP}.
\begin{align}
\label{eqn:P1_QP}
\begin{cases} 
\min \quad  &x^TQx+b^Tx \\
\text{s.t.} & Mx+b \geq 0 \\
&x \geq 0 \\
\end{cases} & \text{, where} \quad
Q= \begin{bmatrix}
A		&0		&0		&0 \\
0		&0		&0		&0 \\
0		&0		&0		&0 \\
0		&0		&0		&E \end{bmatrix}, 
\end{align}
The properties of $A \in \mathbb{R}^{r \times r}$ and $E$ imply that $Q \in \mathbb{R}^{p \times p}$ is diagonal and positive semi-definite. We are interested in the behavior of the solution $\hat{x}$ under perturbation of the diagonal elements of $M$ and $Q$. The diagonal elements $M_{ii}=Q_{ii} \in \mathcal{A}^i \subset \mathbb{R}_0^+$. $\mathcal{A}^i$ is referred to as the admissible set for element $i$, for which the problem remains feasible for all $M_{ii} \in \mathcal{A}^i$. We rewrite Problem \eqref{eqn:P1_QP} in a more general setting as a parametric convex program:
\begin{align} 
\label{eqn:P1_Palpha}
P^\ast(M_{ii})=\begin{cases} \min \limits_{x_i,x_J} &f(x_i,x_J, M_{ii}) \\
\text{s.t.}  &g(x_J) + k(x_i,M_{ii}) \geq 0,\end{cases}
\end{align}
where $J = \{1, \ldots, p\} \setminus \{i\}$. Problem \eqref{eqn:P1_Palpha} can be described in the framework of \eqref{eqn:P1_QP} by setting
\begin{align} 
f(x_i,x_J, M_{ii}) 	&= \sum \limits_{j \in J} \left(Q_{jj} x_j^2 + b_{j} x_j \right) + M_{ii} x_i^2 + b_i x_i\\
g(x_J) + k(x_i,M_{ii}) 	&= 	\begin{bmatrix} M_{ii} & M_{iJ} \\ M_{Ji} & M_{JJ} \end{bmatrix} \cdot \begin{bmatrix} x_i\\x_J \end{bmatrix} + \begin{bmatrix} b_i\\b_J \end{bmatrix} \\
								&= 	\begin{bmatrix} M_{iJ} \\ M_{JJ} \end{bmatrix} x_J + \begin{bmatrix} M_{ii} \\ M_{Ji} \end{bmatrix} x_i +  \begin{bmatrix} b_i\\b_J \end{bmatrix}.
\end{align}
We base our analysis on Problem \eqref{eqn:P1_Palpha} and on the following assumptions:
\begin{enumerate}[(i)]
	\item $x_J \mapsto g(x_J)$ convex,
	\item $x_i \mapsto k(x_i,M_{ii})$ convex for all $M_{ii} \in \mathbb{R}_{M_{ii}\geq 0}$,
	\item $(x_i,x_J) \mapsto f(x_i,x_J,M_{ii})$ convex for all $M_{ii} \in \mathbb{R}_{M_{ii}\geq 0}$,
\end{enumerate}
	
\paragraph{Claim}
If $x_i \mapsto f(x_i,x_J,M_{ii})$ is a strictly convex mapping for all $x_J$, the optimizer of Problem \eqref{eqn:P1_Palpha} admits a unique $i$-component.

\paragraph{Proof}
For the sake of contradiction, suppose $\begin{bmatrix} \hat{x}^1_i \\ \hat{x}^1_J \end{bmatrix}$ and $\begin{bmatrix} \hat{x}^2_i \\ \hat{x}^2_J \end{bmatrix}$ are two optimizers, where $\hat{x}^1_i \neq \hat{x}^2_i$. Then, clearly $\begin{bmatrix}\frac{1}{2}\hat{x}^1_i + \frac{1}{2}\hat{x}^2_i \\ \frac{1}{2}\hat{x}^1_J + \frac{1}{2}\hat{x}^2_J \end{bmatrix}$ is a feasible point with a lower objective value. $\Box$

In the specific case of Problem \eqref{eqn:P1_QP}, note that (i) strict convexity in $x_i \mapsto f(x_i,x_J,M_{ii})$ for all $x_J$ implies $M_{ii}>0$, and (ii) the problem is feasible in particular for all $A_{ii} \in \mathbb{R}_{A_{ii}\geq 0}$, $i \in \{1, \ldots, r\}$, since there always exists a solution $\hat{x}$ with $x_i=0$ for all $i \in \{1, \ldots, r\}$ (Section \ref{sec:P1_ExistenceOfASolution}). 

\section{Available data sources} \label{app:P1_AvailableDataSources}
\begin{table}[htbp]
\centering
\caption{Currently available data for model parameters.}
\label{tab:P1_ListOfData_Params} 
\begin{tabular}{p{0.22\textwidth}p{0.68\textwidth}}
\hline 
\textbf{Parameter}  &  \textbf{Data source} \\
\hline 
$\overline{\CAP}^P$, $\overline{\CAP}^{PT}$ &  Estimated based on data of $s^{P}$ \\ 
$\overline{\CAP}^A$ & ENTSO-G \citep{ENTSO-G2012a}, \citet{Egging2008}, BP \citep{BP2013}, EIA \citep{U.S.EnergyInformationAdministration}, Energy Charter \citep{EnergyCharter} \\ 
$\overline{\CAP}^L$ & GLE \citep{GasLNGEuropeGLE2012}, GIIGNL \citep{InternationalGroupofLiquefiedNaturalGasImportersGIIGNL2012}, BP \citep{BP2013} \\ 
$\overline{\CAP}^R$ & GLE \citep{GasLNGEuropeGLE2012}, GIIGNL \citep{InternationalGroupofLiquefiedNaturalGasImportersGIIGNL2012}, BP \citep{BP2013} \\ 
$\overline{\CAP}^{I}$, $\overline{\CAP}^{IT}$, $\overline{\CAP}^{X}$& GSE \citep{GasStorageEuropeGSE2012}, CGA \citep{FirstEnergyCapitalCorp2014a}, EIA \citep{U.S.EnergyInformationAdministration}, TPAC \citep{TurkishPetroleumCorporation}, Inogate \citep{Inogate},  \citet{Yoshizaki2008}, IEA \citep{InternationalEnergyAgencyIEA2012}\\ 
\hline
$\LINC^P$, $\QUAC^P$ & \citet{Egging2008} \\ 
$\LINC^A$ & IEA \citep{InternationalEnergyAgencyIEA2009}, ILF \citep{Schwimmbeck2008} \\ 
$\LINC^L$ & \citet{Egging2008} \\ 
$\LINC^B$ & \citet{Egging2008} \\ 
$\LINC^R$ & \citet{Egging2008} \\ 
$\LINC^S$ & \citet{Egging2008} \\ 
\hline
$\LOSS^A$ & IGU \citep{InternationalGasUnionIGU2012} \\
$\LOSS^L$ & Petrowiki \citep{PetroWiki2012} \\ 
$\LOSS^B$ & \citet{Egging2008}, Petrowiki \citep{PetroWiki2012}, GIIGNL \citep{InternationalGroupofLiquefiedNaturalGasImportersGIIGNL2012}, \citet{Sea-Distances.org}, \citet{Searates.com} \\ 
$\LOSS^R$ & Petrowiki \citep{PetroWiki2012}, \citet{Egging2008} \\
$\LOSS^S$ & \citet{Egging2008}\\ 
\hline
$\INT^C$, $\SLP^C$ & \citet{Lise2008}, UN \citep{UnitedNations} \\
\hline 
\end{tabular}
\end{table}

\begin{table}[htbp]
\centering
\caption{Currently available data for model variables. This data was used for the model calibration.}
\label{tab:P1_ListOfData_Vars} 
\begin{tabular}{p{0.2\textwidth}p{0.7\textwidth}ll}
\hline
\textbf{Variable} &  \textbf{Data source} \\
\hline 
$q^{C}, \sum \limits_{f \in \mathcal{F}(n)} q^{C}_{fnt}, $ & OECD \citep{OECD}, BP \citep{BP2013}, EC \citep{EuropeanCommission}, IEA \citep{InternationalEnergyAgencyIEA2013}, EIA \citep{U.S.EnergyInformationAdministration}, UN \citep{UnitedNations}\\
$\lambda^{C}$ & IGU \citep{InternationalGasUnionIGU2013}, EC \citep{EuropeanCommission}, IEA \citep{InternationalEnergyAgencyIEA2013}, EIA \citep{U.S.EnergyInformationAdministration} \\ 
\hline 
\end{tabular}
\end{table}

\clearpage
\bibliography{Bibliography}

\begin{thebibliography}{51}
\expandafter\ifx\csname natexlab\endcsname\relax\def\natexlab#1{#1}\fi
\providecommand{\url}[1]{\texttt{#1}}
\providecommand{\href}[2]{#2}
\providecommand{\path}[1]{#1}
\providecommand{\DOIprefix}{doi:}
\providecommand{\ArXivprefix}{arXiv:}
\providecommand{\URLprefix}{URL: }
\providecommand{\Pubmedprefix}{pmid:}
\providecommand{\doi}[1]{\href{http://dx.doi.org/#1}{\path{#1}}}
\providecommand{\Pubmed}[1]{\href{pmid:#1}{\path{#1}}}
\providecommand{\bibinfo}[2]{#2}
\ifx\xfnm\relax \def\xfnm[#1]{\unskip,\space#1}\fi
\bibitem[{Abolhassani et~al.(2014)Abolhassani, Bateni, Hajiaghayi, Mahini \&
  Sawant}]{Abolhassani2014}
\bibinfo{author}{Abolhassani, M.}, \bibinfo{author}{Bateni, M.},
  \bibinfo{author}{Hajiaghayi, M.}, \bibinfo{author}{Mahini, H.}, \&
  \bibinfo{author}{Sawant, A.} \bibinfo{year}{2014}.
\newblock \bibinfo{title}{{Network Cournot Competition}}.
\newblock In \bibinfo{editor}{T.-W. Liu}, \bibinfo{editor}{Q.~Qi}, \&
  \bibinfo{editor}{Y.~Ye} (Eds.), {\it \bibinfo{booktitle}{Web and Internet
  Economics}\/} (pp. \bibinfo{pages}{15--29}).
\newblock \bibinfo{address}{Cham}: \bibinfo{publisher}{Springer}.
\newblock \href{http://arxiv.org/abs/arXiv:1405.1794v1}{\tt
  arXiv:arXiv:1405.1794v1}.
\bibitem[{BP(2013)}]{BP2013}
\bibinfo{author}{BP} \bibinfo{year}{2013}.
\newblock {\it \bibinfo{title}{{BP statistical review of world energy June
  2013}}\/}.
\newblock \bibinfo{type}{Technical Report} \bibinfo{address}{retrieved from
  \url{http://www.bp.com/en/global/corporate/about-bp/energy-economics/statistical-review-of-world-energy/2013-in-review.html}}.
\bibitem[{Carvalho et~al.(2014)Carvalho, Buzna, Bono, Masera, Arrowsmith \&
  Helbing}]{Carvalho2014}
\bibinfo{author}{Carvalho, R.}, \bibinfo{author}{Buzna, L.},
  \bibinfo{author}{Bono, F.}, \bibinfo{author}{Masera, M.},
  \bibinfo{author}{Arrowsmith, D.~K.}, \& \bibinfo{author}{Helbing, D.}
  \bibinfo{year}{2014}.
\newblock \bibinfo{title}{{Resilience of natural gas networks during conflicts,
  crises and disruptions.}}
\newblock {\it \bibinfo{journal}{PloS one}\/},  {\it \bibinfo{volume}{9}\/},
  \bibinfo{pages}{e90265}. \DOIprefix\doi{10.1371/journal.pone.0090265}.
\bibitem[{Chyong \& Hobbs(2014)}]{Chyong2014}
\bibinfo{author}{Chyong, C.~K.}, \& \bibinfo{author}{Hobbs, B.~F.}
  \bibinfo{year}{2014}.
\newblock \bibinfo{title}{{Strategic Eurasian natural gas market model for
  energy security and policy analysis: Formulation and application to South
  Stream}}.
\newblock {\it \bibinfo{journal}{Energy Economics}\/},  {\it
  \bibinfo{volume}{44}\/}, \bibinfo{pages}{198--211}.
  \DOIprefix\doi{10.1016/j.eneco.2014.04.006}.
\bibitem[{Cottle et~al.(1992)Cottle, Pang \& Stone}]{Cottle1992}
\bibinfo{author}{Cottle, R.~W.}, \bibinfo{author}{Pang, J.-S.}, \&
  \bibinfo{author}{Stone, R.~E.} \bibinfo{year}{1992}.
\newblock {\it \bibinfo{title}{{The linear complementarity problem}}\/}.
\newblock (\bibinfo{edition}{{SIAM}} ed.).
\newblock \bibinfo{address}{Philadelphia}: \bibinfo{publisher}{Academic Press}.
\bibitem[{D\'{\i}az et~al.(2010)D\'{\i}az, Villar, Campos \&
  Reneses}]{Diaz2010}
\bibinfo{author}{D\'{\i}az, C.~a.}, \bibinfo{author}{Villar, J.},
  \bibinfo{author}{Campos, F.~A.}, \& \bibinfo{author}{Reneses, J.}
  \bibinfo{year}{2010}.
\newblock \bibinfo{title}{{Electricity market equilibrium based on conjectural
  variations}}.
\newblock {\it \bibinfo{journal}{Electric Power Systems Research}\/},  {\it
  \bibinfo{volume}{80}\/}, \bibinfo{pages}{1572--1579}.
  \DOIprefix\doi{10.1016/j.epsr.2010.07.012}.
\bibitem[{Dieckh\"{o}ner et~al.(2013)Dieckh\"{o}ner, Lochner \&
  Lindenberger}]{Dieckhoner2013}
\bibinfo{author}{Dieckh\"{o}ner, C.}, \bibinfo{author}{Lochner, S.}, \&
  \bibinfo{author}{Lindenberger, D.} \bibinfo{year}{2013}.
\newblock \bibinfo{title}{{European natural gas infrastructure: The impact of
  market developments on gas flows and physical market integration}}.
\newblock {\it \bibinfo{journal}{Applied Energy}\/},  {\it
  \bibinfo{volume}{102}\/}, \bibinfo{pages}{994--1003}.
  \DOIprefix\doi{10.1016/j.apenergy.2012.06.021}.
\bibitem[{Egging(2013)}]{Egging2013}
\bibinfo{author}{Egging, R.} \bibinfo{year}{2013}.
\newblock \bibinfo{title}{{Benders Decomposition for multi-stage stochastic
  mixed complementarity problems -- Applied to a global natural gas market
  model}}.
\newblock {\it \bibinfo{journal}{European Journal of Operational Research}\/},
  {\it \bibinfo{volume}{226}\/}, \bibinfo{pages}{341--353}.
  \DOIprefix\doi{10.1016/j.ejor.2012.11.024}.
\bibitem[{Egging et~al.(2008)Egging, Gabriel, Holz \& Zhuang}]{Egging2008}
\bibinfo{author}{Egging, R.}, \bibinfo{author}{Gabriel, S.~A.},
  \bibinfo{author}{Holz, F.}, \& \bibinfo{author}{Zhuang, J.}
  \bibinfo{year}{2008}.
\newblock \bibinfo{title}{{A complementarity model for the European natural gas
  market}}.
\newblock {\it \bibinfo{journal}{Energy Policy}\/},  {\it
  \bibinfo{volume}{36}\/}, \bibinfo{pages}{2385--2414}.
  \DOIprefix\doi{10.1016/j.enpol.2008.01.044}.
\bibitem[{Egging et~al.(2010)Egging, Holz \& Gabriel}]{Egging2010}
\bibinfo{author}{Egging, R.}, \bibinfo{author}{Holz, F.}, \&
  \bibinfo{author}{Gabriel, S.~A.} \bibinfo{year}{2010}.
\newblock \bibinfo{title}{{The World Gas Model: A multi-period mixed
  complementarity model for the global natural gas market}}.
\newblock {\it \bibinfo{journal}{Energy}\/},  {\it \bibinfo{volume}{35}\/},
  \bibinfo{pages}{4016--4029}. \DOIprefix\doi{10.1016/j.energy.2010.03.053}.
\bibitem[{Egging(2010)}]{Egging2010a}
\bibinfo{author}{Egging, R.~G.} \bibinfo{year}{2010}.
\newblock {\it \bibinfo{title}{{Multi-period natural gas market modeling:
  Applications, stochastic extensions and solution approaches}}\/}.
\newblock Ph.D. thesis University of Maryland, College Park.
\bibitem[{{Energy Charter Secretariat}(2009)}]{EnergyCharter}
\bibinfo{author}{{Energy Charter Secretariat}} \bibinfo{year}{2009}.
\newblock \bibinfo{title}{{Master plan: Ukrainian gas transmission system
  (UGTS), priority objects, modernisation and reconstruction}}.
\bibitem[{ENTSO-G(2012)}]{ENTSO-G2012a}
\bibinfo{author}{ENTSO-G} \bibinfo{year}{2012}.
\newblock {\it \bibinfo{title}{{The European natural gas network: Capacities at
  cross-border points on the primary market}}\/}.
\newblock \bibinfo{type}{Technical Report} \bibinfo{address}{retrieved from
  \url{http://www.entsog.eu/maps/transmission-capacity-map/2012}}.
\bibitem[{{European Commission}(2015)}]{EuropeanCommission}
\bibinfo{author}{{European Commission}} \bibinfo{year}{2015}.
\newblock \bibinfo{title}{{Eurostat}}.
\newblock
  \bibinfo{howpublished}{\url{http://ec.europa.eu/eurostat/de/data/database}}.
\newblock \bibinfo{annote}{{[Online; accessed 13-July-2015]}}.
\bibitem[{{First Energy Capital Corp}(2014)}]{FirstEnergyCapitalCorp2014a}
\bibinfo{author}{{First Energy Capital Corp}} \bibinfo{year}{2014}.
\newblock \bibinfo{title}{{Natural gas storage -- Canada}}.
\newblock
  \bibinfo{howpublished}{\url{http://www.cga.ca/wp-content/uploads/2011/02/Chart-1-Natural-Gas-Storage38.pdf}}.
\newblock \bibinfo{annote}{{[Online; accessed 09-October-2014]}}.
\bibitem[{Gabriel et~al.(2005)Gabriel, Kiet \& Zhuang}]{Gabriel2005a}
\bibinfo{author}{Gabriel, S.~a.}, \bibinfo{author}{Kiet, S.}, \&
  \bibinfo{author}{Zhuang, J.} \bibinfo{year}{2005}.
\newblock \bibinfo{title}{{A mixed complementarity-based equilibrium model of
  natural gas markets}}.
\newblock {\it \bibinfo{journal}{Operations Research}\/},  {\it
  \bibinfo{volume}{53}\/}, \bibinfo{pages}{799--818}.
  \DOIprefix\doi{10.1287/opre.1040.0199}.
\bibitem[{Gabriel et~al.(2012)Gabriel, Rosendahl, Egging, Avetisyan \&
  Siddiqui}]{Gabriel2012}
\bibinfo{author}{Gabriel, S.~A.}, \bibinfo{author}{Rosendahl, K.~E.},
  \bibinfo{author}{Egging, R.}, \bibinfo{author}{Avetisyan, H.~G.}, \&
  \bibinfo{author}{Siddiqui, S.} \bibinfo{year}{2012}.
\newblock \bibinfo{title}{{Cartelization in gas markets: Studying the potential
  for a ``Gas OPEC"}}.
\newblock {\it \bibinfo{journal}{Energy Economics}\/},  {\it
  \bibinfo{volume}{34}\/}, \bibinfo{pages}{137--152}.
  \DOIprefix\doi{10.1016/j.eneco.2011.05.014}.
\bibitem[{Gabriel et~al.(2009)Gabriel, Zhuang \& Egging}]{Gabriel2009}
\bibinfo{author}{Gabriel, S.~A.}, \bibinfo{author}{Zhuang, J.}, \&
  \bibinfo{author}{Egging, R.} \bibinfo{year}{2009}.
\newblock \bibinfo{title}{{Solving stochastic complementarity problems in
  energy market modeling using scenario reduction}}.
\newblock {\it \bibinfo{journal}{European Journal of Operational Research}\/},
  {\it \bibinfo{volume}{197}\/}, \bibinfo{pages}{1028--1040}.
  \DOIprefix\doi{10.1016/j.ejor.2007.12.046}.
\bibitem[{{Gas LNG Europe (gle)}(2012)}]{GasLNGEuropeGLE2012}
\bibinfo{author}{{Gas LNG Europe (gle)}} \bibinfo{year}{2012}.
\newblock {\it \bibinfo{title}{{LNG map}}\/}.
\newblock \bibinfo{type}{Technical Report} \bibinfo{address}{retrieved from
  \url{http://www.gie.eu/index.php/maps-data/lng-map}}.
\bibitem[{{Gas Storage Europe (gse)}(2012)}]{GasStorageEuropeGSE2012}
\bibinfo{author}{{Gas Storage Europe (gse)}} \bibinfo{year}{2012}.
\newblock {\it \bibinfo{title}{{Storage Map}}\/}.
\newblock \bibinfo{type}{Technical Report} \bibinfo{address}{retrieved from
  \url{http://www.gie.eu/index.php/maps-data/gse-storage-map}}.
\bibitem[{Huppmann(2013)}]{Huppmann2013a}
\bibinfo{author}{Huppmann, D.} \bibinfo{year}{2013}.
\newblock \bibinfo{title}{{Endogenous production capacity investment in natural
  gas market equilibrium models}}.
\newblock {\it \bibinfo{journal}{European Journal of Operational Research}\/},
  {\it \bibinfo{volume}{231}\/}, \bibinfo{pages}{503--506}.
  \DOIprefix\doi{10.1016/j.ejor.2013.05.048}.
\bibitem[{Huppmann \& Egging(2014)}]{Huppmann2014}
\bibinfo{author}{Huppmann, D.}, \& \bibinfo{author}{Egging, R.}
  \bibinfo{year}{2014}.
\newblock \bibinfo{title}{{Market power, fuel substitution and infrastructure
  – A large-scale equilibrium model of global energy markets}}.
\newblock {\it \bibinfo{journal}{Energy}\/},  {\it \bibinfo{volume}{75}\/},
  \bibinfo{pages}{483--500}. \DOIprefix\doi{10.1016/j.energy.2014.08.004}.
\bibitem[{IEA(2009)}]{InternationalEnergyAgencyIEA2009}
\bibinfo{author}{IEA} \bibinfo{year}{2009}.
\newblock {\it \bibinfo{title}{{World energy outlook 2009}}\/}.
\newblock \bibinfo{type}{Technical Report} \bibinfo{address}{retrieved from
  \url{http://www.worldenergyoutlook.org/media/weowebsite/2009/WEO2009.pdf}}.
\bibitem[{IEA(2012)}]{InternationalEnergyAgencyIEA2012}
\bibinfo{author}{IEA} \bibinfo{year}{2012}.
\newblock {\it \bibinfo{title}{{Gas pricing and regulation: China's challenges
  and IEA experience}}\/}.
\newblock \bibinfo{type}{Technical Report} \bibinfo{address}{retrieved from
  \url{http://www.iea.org/publications/freepublications/publication/chinagasreport_final_web.pdf}}.
\bibitem[{IEA(2013)}]{InternationalEnergyAgencyIEA2013}
\bibinfo{author}{IEA} \bibinfo{year}{2013}.
\newblock {\it \bibinfo{title}{{IEA statistics: Natural gas information}}\/}.
\newblock \bibinfo{type}{Technical Report}.
\newblock \DOIprefix\doi{10.1787/nat_gas-2013-en}.
\bibitem[{Inogate()}]{Inogate}
\bibinfo{author}{Inogate}  (n.d.).
\newblock \bibinfo{title}{{Inogate umbrella agreement: Republic of Belarus}}.
\newblock
  \bibinfo{howpublished}{\url{http://www2.inogate.org/html/countries/belarus.htm}}.
\newblock \bibinfo{annote}{{[Online; accessed 03-December-2014]}}.
\bibitem[{{International Gas Union}(2012)}]{InternationalGasUnionIGU2012}
\bibinfo{author}{{International Gas Union}} \bibinfo{year}{2012}.
\newblock \bibinfo{title}{{Natural gas: Facts \& figures}}.
\bibitem[{{International Gas Union}(2013)}]{InternationalGasUnionIGU2013}
\bibinfo{author}{{International Gas Union}} \bibinfo{year}{2013}.
\newblock {\it \bibinfo{title}{{Wholesale gas price survey -- 2013 edition: A
  global review of price formation mechanisms 2005 -2012}}\/}.
\newblock \bibinfo{type}{Technical Report} \bibinfo{address}{retrieved from
  \url{http://members.igu.org/news/igu-launces-the-wholesale-gas-price-survey-2013-edition}}.
\bibitem[{{International Group of Liquefied Natural Gas Importers
  (GIIGNL)}(2012)}]{InternationalGroupofLiquefiedNaturalGasImportersGIIGNL2012}
\bibinfo{author}{{International Group of Liquefied Natural Gas Importers
  (GIIGNL)}} \bibinfo{year}{2012}.
\newblock {\it \bibinfo{title}{{The LNG Industry}}\/}.
\newblock \bibinfo{type}{Technical Report} \bibinfo{address}{retrieved from
  \url{http://www.giignl.org/system/files/publication/giignl_the_lng_industry_2012.pdf}}.
\bibitem[{de~Joode \& \"{O}zdemir(2010)}]{DeJoode2010}
\bibinfo{author}{de~Joode, J.}, \& \bibinfo{author}{\"{O}zdemir, O.}
  \bibinfo{year}{2010}.
\newblock \bibinfo{title}{{Demand for seasonal gas storage in northwest Europe
  until 2030: Simulation results with a dynamic model}}.
\newblock {\it \bibinfo{journal}{Energy Policy}\/},  {\it
  \bibinfo{volume}{38}\/}, \bibinfo{pages}{5817--5829}.
  \DOIprefix\doi{10.1016/j.enpol.2010.05.032}.
\bibitem[{Kamiński(2011)}]{Kaminski2011}
\bibinfo{author}{Kamiński, J.} \bibinfo{year}{2011}.
\newblock \bibinfo{title}{{Market power in a coal-based power generation
  sector: The case of Poland}}.
\newblock {\it \bibinfo{journal}{Energy}\/},  {\it \bibinfo{volume}{36}\/},
  \bibinfo{pages}{6634--6644}. \DOIprefix\doi{10.1016/j.energy.2011.08.048}.
\bibitem[{Lise \& Hobbs(2008)}]{Lise2008b}
\bibinfo{author}{Lise, W.}, \& \bibinfo{author}{Hobbs, B.~F.}
  \bibinfo{year}{2008}.
\newblock \bibinfo{title}{{Future evolution of the liberalised European gas
  market: Simulation results with a dynamic model}}.
\newblock {\it \bibinfo{journal}{Energy}\/},  {\it \bibinfo{volume}{33}\/},
  \bibinfo{pages}{989--1004}. \DOIprefix\doi{10.1016/j.energy.2008.02.012}.
\bibitem[{Lise et~al.(2008)Lise, Hobbs \& van Oostvoorn}]{Lise2008}
\bibinfo{author}{Lise, W.}, \bibinfo{author}{Hobbs, B.~F.}, \&
  \bibinfo{author}{van Oostvoorn, F.} \bibinfo{year}{2008}.
\newblock \bibinfo{title}{{Natural gas corridors between the EU and its main
  suppliers: Simulation results with the dynamic GASTALE model}}.
\newblock {\it \bibinfo{journal}{Energy Policy}\/},  {\it
  \bibinfo{volume}{36}\/}, \bibinfo{pages}{1890--1906}. \DOIprefix\doi{DOI
  10.1016/j.enpol.2008.01.042}.
\bibitem[{Lochner(2011{\natexlab{a}})}]{Lochner2011}
\bibinfo{author}{Lochner, S.} \bibinfo{year}{2011}{\natexlab{a}}.
\newblock \bibinfo{title}{{Identification of congestion and valuation of
  transport infrastructures in the European natural gas market}}.
\newblock {\it \bibinfo{journal}{Energy}\/},  {\it \bibinfo{volume}{36}\/},
  \bibinfo{pages}{2483--2492}. \DOIprefix\doi{10.1016/j.energy.2011.01.040}.
\bibitem[{Lochner(2011{\natexlab{b}})}]{Lochner2011a}
\bibinfo{author}{Lochner, S.} \bibinfo{year}{2011}{\natexlab{b}}.
\newblock \bibinfo{title}{{Modeling the European natural gas market during the
  2009 Russian-Ukrainian gas conflict: Ex-post simulation and analysis}}.
\newblock {\it \bibinfo{journal}{Journal of Natural Gas Science and
  Engineering}\/},  {\it \bibinfo{volume}{3}\/}, \bibinfo{pages}{341--348}.
  \DOIprefix\doi{10.1016/j.jngse.2011.01.003}.
\bibitem[{Lochner \& Dieckh\"{o}ner(2012)}]{Lochner2012}
\bibinfo{author}{Lochner, S.}, \& \bibinfo{author}{Dieckh\"{o}ner, C.}
  \bibinfo{year}{2012}.
\newblock \bibinfo{title}{{Civil unrest in North Africa\textendash Risks for
  natural gas supply?}}
\newblock {\it \bibinfo{journal}{Energy Policy}\/},  {\it
  \bibinfo{volume}{45}\/}, \bibinfo{pages}{167--175}.
  \DOIprefix\doi{10.1016/j.enpol.2012.02.009}.
\bibitem[{Moryadee et~al.(2014{\natexlab{a}})Moryadee, Gabriel \&
  Avetisyan}]{Moryadee2014}
\bibinfo{author}{Moryadee, S.}, \bibinfo{author}{Gabriel, S.~A.}, \&
  \bibinfo{author}{Avetisyan, H.~G.} \bibinfo{year}{2014}{\natexlab{a}}.
\newblock \bibinfo{title}{{Investigating the potential effects of U.S. LNG
  exports on global natural gas markets}}.
\newblock {\it \bibinfo{journal}{Energy Strategy Reviews}\/},  {\it
  \bibinfo{volume}{2}\/}, \bibinfo{pages}{273--288}.
  \DOIprefix\doi{10.1016/j.esr.2013.12.004}.
\bibitem[{Moryadee et~al.(2014{\natexlab{b}})Moryadee, Gabriel \&
  Rehulka}]{Moryadee2014a}
\bibinfo{author}{Moryadee, S.}, \bibinfo{author}{Gabriel, S.~a.}, \&
  \bibinfo{author}{Rehulka, F.} \bibinfo{year}{2014}{\natexlab{b}}.
\newblock \bibinfo{title}{{The influence of the Panama Canal on global gas
  trade}}.
\newblock {\it \bibinfo{journal}{Journal of Natural Gas Science and
  Engineering}\/},  {\it \bibinfo{volume}{20}\/}, \bibinfo{pages}{161--174}.
  \DOIprefix\doi{10.1016/j.jngse.2014.06.015}.
\bibitem[{Neuhoff et~al.(2005)Neuhoff, Barquin, Boots, Ehrenmann, Hobbs,
  Rijkers \& V\'{a}zquez}]{Neuhoff2005}
\bibinfo{author}{Neuhoff, K.}, \bibinfo{author}{Barquin, J.},
  \bibinfo{author}{Boots, M.~G.}, \bibinfo{author}{Ehrenmann, A.},
  \bibinfo{author}{Hobbs, B.~F.}, \bibinfo{author}{Rijkers, F. a.~M.}, \&
  \bibinfo{author}{V\'{a}zquez, M.} \bibinfo{year}{2005}.
\newblock \bibinfo{title}{{Network-constrained Cournot models of liberalized
  electricity markets: The devil is in the details}}.
\newblock {\it \bibinfo{journal}{Energy Economics}\/},  {\it
  \bibinfo{volume}{27}\/}, \bibinfo{pages}{495--525}.
  \DOIprefix\doi{10.1016/j.eneco.2004.12.001}.
\bibitem[{OECD()}]{OECD}
\bibinfo{author}{OECD}  (n.d.).
\newblock \bibinfo{title}{{iLibrary}}.
\newblock
  \bibinfo{howpublished}{\url{http://www.oecd-ilibrary.org/statistics}}.
\newblock \bibinfo{annote}{{[Online; accessed 09-October-2014]}}.
\bibitem[{Perloff et~al.(2007)Perloff, Karp \& Golan}]{Perloff2007}
\bibinfo{author}{Perloff, J.~M.}, \bibinfo{author}{Karp, L.~S.}, \&
  \bibinfo{author}{Golan, A.} \bibinfo{year}{2007}.
\newblock \bibinfo{title}{{Industry Models of Market Power}}.
\newblock In {\it \bibinfo{booktitle}{Estimating Market Power and
  Strategies}\/} (pp. \bibinfo{pages}{42--73}).
\newblock \bibinfo{publisher}{Cambridge University Press}.
\newblock \DOIprefix\doi{http://dx.doi.org/10.1017/CBO9780511753985.004}.
\bibitem[{PetroWiki(2012)}]{PetroWiki2012}
\bibinfo{author}{PetroWiki} \bibinfo{year}{2012}.
\newblock \bibinfo{title}{{Efficiency losses in the LNG value chain}}.
\newblock \bibinfo{annote}{{[Online; accessed 08-December-2014]}}.
\bibitem[{Schwimmbeck(2008)}]{Schwimmbeck2008}
\bibinfo{author}{Schwimmbeck, R.~G.} \bibinfo{year}{2008}.
\newblock \bibinfo{title}{{LNG and pipeline}}.
\newblock In {\it \bibinfo{booktitle}{3rd Pipeline Technology Conference}\/}.
\newblock \bibinfo{address}{Retrieved from
  \url{http://www.pipeline-conference.com/download/1158/ptc_2008_schwimmbeck.pdf}}.
\bibitem[{Sea-distances.org()}]{Sea-Distances.org}
\bibinfo{author}{Sea-distances.org}  (n.d.).
\newblock \bibinfo{title}{{Port distances}}.
\newblock \bibinfo{howpublished}{\url{http://www.sea-distances.org/}}.
\newblock \bibinfo{annote}{{[Online; accessed 03-December-2014]}}.
\bibitem[{SeaRates.com()}]{Searates.com}
\bibinfo{author}{SeaRates.com}  (n.d.).
\newblock \bibinfo{title}{{Distances and time}}.
\newblock
  \bibinfo{howpublished}{\url{http://www.searates.com/reference/portdistance/}}.
\newblock \bibinfo{annote}{{[Online; accessed 03-December-2014]}}.
\bibitem[{Tremblay \& Tremblay(2012)}]{Tremblay2012}
\bibinfo{author}{Tremblay, V.~J.}, \& \bibinfo{author}{Tremblay, C.~H.}
  \bibinfo{year}{2012}.
\newblock {\it \bibinfo{title}{{New Perspectives on Industrial Organization:
  With Contributions from Behavioral Economics and Game Theory}}\/}.
\newblock \bibinfo{address}{New York}: \bibinfo{publisher}{Springer}.
\newblock \DOIprefix\doi{10.1007/978-1-4614-3241-8}.
\bibitem[{{Turkish Petroleum Corporation}()}]{TurkishPetroleumCorporation}
\bibinfo{author}{{Turkish Petroleum Corporation}}  (n.d.).
\newblock \bibinfo{title}{{Natural gas storage}}.
\newblock \bibinfo{howpublished}{\url{http://www.tpao.gov.tr/eng/?tp=m&id=84}}.
\newblock \bibinfo{annote}{{[Online; accessed 03-December-2014]}}.
\bibitem[{{United Nations}()}]{UnitedNations}
\bibinfo{author}{{United Nations}}  (n.d.).
\newblock \bibinfo{title}{{UN data}}.
\newblock \bibinfo{howpublished}{\url{http://data.un.org/}}.
\newblock \bibinfo{annote}{{[Online; accessed 09-October-2014]}}.
\bibitem[{{U.S. Energy Information Administration
  (EIA)}()}]{U.S.EnergyInformationAdministration}
\bibinfo{author}{{U.S. Energy Information Administration (EIA)}}  (n.d.).
\newblock \bibinfo{title}{{Natural gas data}}.
\newblock \bibinfo{howpublished}{\url{http://www.eia.gov/naturalgas/data.cfm}.}
\newblock \bibinfo{annote}{{[Online; accessed 09-October-2014]}}.
\bibitem[{Yoshizaki et~al.(2008)Yoshizaki, Sato, Fukagawa, Sugiyama, Takagi \&
  Jono}]{Yoshizaki2008}
\bibinfo{author}{Yoshizaki, K.}, \bibinfo{author}{Sato, N.},
  \bibinfo{author}{Fukagawa, H.}, \bibinfo{author}{Sugiyama, H.},
  \bibinfo{author}{Takagi, G.}, \& \bibinfo{author}{Jono, T.}
  \bibinfo{year}{2008}.
\newblock \bibinfo{title}{{Utilization of underground gas storage (UGS) in
  Japan}}.
\bibitem[{Zhuang \& Gabriel(2008)}]{Zhuang2008}
\bibinfo{author}{Zhuang, J.}, \& \bibinfo{author}{Gabriel, S.~A.}
  \bibinfo{year}{2008}.
\newblock \bibinfo{title}{{A complementarity model for solving stochastic
  natural gas market equilibria}}.
\newblock {\it \bibinfo{journal}{Energy Economics}\/},  {\it
  \bibinfo{volume}{30}\/}, \bibinfo{pages}{113--147}.
  \DOIprefix\doi{10.1016/j.eneco.2006.09.004}.

\end{thebibliography}
\end{document}